\input amstex
\documentstyle{amsppt}
\magnification=\magstep1

\pageheight{9.0truein}
\pagewidth{6.5truein}

\hyphenation{uni-ser-ial}

\input xy
\xyoption{matrix}\xyoption{arrow}\xyoption{curve}

\long\def\ignore#1{#1}

\def\uloopr#1{\ar@'{@+{[0,0]+(-4,5)} @+{[0,0]+(0,10)} @+{[0,0]+(4,5)}}
  ^{#1}}
\def\dloopr#1{\ar@'{@+{[0,0]+(-4,-5)} @+{[0,0]+(0,-10)} @+{[0,0]+(4,-5)}}
  _{#1}}
\def\rloopd#1{\ar@'{@+{[0,0]+(5,4)} @+{[0,0]+(10,0)} @+{[0,0]+(5,-4)}}
  ^{#1}}
\def\lloopd#1{\ar@'{@+{[0,0]+(-5,4)} @+{[0,0]+(-10,0)} @+{[0,0]+(-5,-4)}}
  _{#1}}
\def\udotloopr#1{\ar@{{}.>} @'{@+{[0,0]+(-4,5)} @+{[0,0]+(0,10)}
  @+{[0,0]+(4,5)}} ^{#1}}

\def\seq{\mathrel{\widehat{=}}}
\def\la{{\Lambda}}
\def\lamod{\Lambda\text{-}\roman{mod}}
\def \len{\operatorname{length}} 
\def\detour{\mathrel{\wr\wr}}

\def\P{{\Cal P}}

\def\T{{\Cal T}}
\def\S{{\Cal S}}
\def\O{{\Cal O}}
\def\A{{\Cal A}}
\def\B{{\Cal B}}
\def\F{{\Cal F}}
\def\AA{{\Bbb A}}
\def\SS{{\Bbb S}}

\def\soc{\operatorname{soc}}
\def\top{\operatorname{top}}
\def\Ext{\operatorname{Ext^1_\la}}

\def\lfd{\operatorname{l\,fin\,dim}}
\def\rank{\operatorname{rank}}

\topmatter
\title The geometry of uniserial representations of finite dimensional
algebras III: Finite uniserial type
\endtitle
\rightheadtext{Geometry of uniserial representations: Finite uniserial type}
\author Birge Zimmermann Huisgen\endauthor
\address Department of Mathematics, University of California, Santa Barbara,
CA 93106\endaddress
\email birge\@math.ucsb.edu\endemail
\thanks This research was partially supported by a National Science
Foundation grant.\endthanks


\abstract A description is given of those sequences $\SS=
(S(0),S(1),\dots,S(l))$ of simple modules over a finite dimensional algebra
for which there are only finitely many uniserial modules with consecutive
composition factors
$S(0),\dots,S(l)$. Necessary and sufficient conditions for an algebra to
permit only a finite number of isomorphism types of uniserial modules are
derived. The main tools in this investigation are the affine algebraic
varieties parametrizing the uniserial modules
with composition series $\SS$.
\endabstract

\endtopmatter
\document

\head {1. Introduction and setting} \endhead

The purpose of this paper is to develop necessary and sufficient conditions
for a finite dimensional algebra to have `finite uniserial type', that is,
to have only a finite number of isomorphism classes of
uniserial modules. This problem has been raised on many occasions by M.
Auslander since the mid-seventies.  One of the reasons for Auslander's
interest in this condition lies in the fact that for any algebra $\la$ of
infinite uniserial type, Smal\o's 'Brauer-Thrall 1.5' [9] readily
yields the conclusion of the second Brauer-Thrall Conjecture:  Indeed,
infinite uniserial type provides us with an infinite family of non-isomorphic
indecomposable
$\la$-modules of some dimension $d \le \dim_K \la$, where $K$ is the base
field of $\la$, and hence Smal\o's result supplies us with an infinite
sequence of distinct integers, each of which occurs as the $K$-dimension of
infinitely many non-isomorphic indecomposable $\la$-modules. A second motive
-- for both Auslander and the author -- for propagating uniserial
representation theory is as follows.  Being in
some sense the simplest cyclic modules beyond the indecomposable projectives,
the uniserials play a key role as building blocks of the objects in other
classifiable families, and in approaching classes of less accessible modules,
such classifiable families provide a useful basis for comparison.  A concrete
instance of this philosophy can be found in [2], namely in the study of
minimal right approximations of the simple modules by modules of finite
projective dimension (here we use the term `minimal right approximation' in
the sense of Auslander and Reiten [1]).  The present paper is a contribution
to a more general approximation theory which uses uniserial modules as a
starting point. 

Our investigation of split basic algebras $\la$ of finite uniserial type will,
on one hand, provide criteria in terms of quiver and relations and, on the
other hand, zero in on the uniserial representation theory of such algebras,
e.g., describe the sequences of simple modules which arise as sequences of
consecutive composition factors of uniserial left $\la$-modules over such
algebras.  Roughly speaking, it will turn out that, whenever a simple module
occurs with multiplicity larger than $1$ in a given uniserial module, a major
segment of the sequence of composition factors must repeat.  As a
consequence -- even more roughly -- the intricacy of the uniserial
representation theory of an algebra $\la$ of finite uniserial type is
directly proportional to the intricacy of patterns of nested oriented cycles
in the quiver of the algebra.  One of the main points of the paper is to
relate the sequences
$\SS$ of simple $\la$-modules which occur as sequences of composition factors
of only finitely many isomorphism types of uniserials to the geometry of the
affine varietes $V_{\SS}$; in [6], these varieties have been shown to draw a
concise picture of the uniserial modules with composition factor sequence
$\SS$.  We sketch some of the ideas of this prior work in order to outline
the results of the present paper.  Of course, the interest of a
uniserial module viewed by itself is as scant as that of a point on a
curve when considered outside the context of the curve.  In classifying
families of uniserial
$\la$-modules, the interest lies in the number of and interplay among the
free parameters which offer themselves for the description of these modules. 
As was shown in [6], the pertinent tuples of parameters trace the points of
certain affine algebraic varieties, the geometry of which is intimately
related to the relative behavior of the corresponding uniserial modules.  A
bit more precisely: Given any sequence $\SS$ of simple left $\la$-modules,
there exists a finite family $V_{\SS}$ of irreducible algebraic varieties,
together with a canonical map from $V_{\SS}$ onto the set of isomorphism types
of uniserial left $\la$-modules with sequence $\SS$ of consecutive
composition factors which provides effective information on these isomorphism
types.  These varieties depend solely on $\SS$ and the isomorphism
class of $\la$, up to birational equivalence, often even up to isomorphism.  

In particular, we will see in Section 3 that, for an algebra of finite
uniserial type, the varieties
$V_{\SS}$ may be replaced by a single one (not necessarily irreducible) which
is unique up to isomorphism.  This is entailed by the following strong
necessary condition (N) for finite uniserial type over an infinite
field
$K$, which will be established in that section:  Whenever $S$ and
$S'$ are simple left $\la$-modules with $\Ext (S,S') \ne 0$, there exists a
unique uniserial module
$W$ of length $2$ with top $S$ and socle $S'$, and each such uniserial module
$U$ of any length $l \ge 2$ has the property that either $U/J^2 U \cong W$ or
else $J^{l-2} U \cong W$. (Here $J$ denotes the Jacobson radical of $\la$.) 
Given a coordinatization of $\la$, namely a presentation in the form $\la =
K\Gamma / I$, where $\Gamma$ is the quiver of $\la$ and $I$ an admissible
ideal in the path algebra $K\Gamma$,  the validity of Condition (N) can be
translated into conditions on $\Gamma$ and the relations of $\la$ as
follows:  A path $p$ of length $l$  in $K\Gamma$ is called a
{\it mast} of a uniserial left $\la$-module $U$ of composition length $l+1$ if
$pU \ne 0$.  Then (N) says that, whenever a path $p: e \rightarrow e'$ in
$K\Gamma$ of positive length is a mast of a uniserial left $\la$-module and
$\alpha : e \rightarrow e'$ is an arrow, the path $p$ is either equal to
$\alpha$ preceded by a cycle based at $e$, or equal to $\alpha$ followed
by a cycle based at $e'$.  Since, in view of [6], the question of
whether a path
$p$ occurs as a mast can be decided algorithmically, this latter rendering of
condition (N) can be readily checked in terms of quiver and relations. A host
of consequences (see Theorem 3.3(II) and Corollary 3.4) shows
how strongly condition (N) impinges on the quiver
$\Gamma$, as well as on the sequences of simple modules occurring as
sequences of consecutive composition factors of uniserial modules.  For
instance, if $S$ is a simple left $\la$-module with
$\Ext(S,S) \ne 0$, then $S$ can only occur in a `connected segment' within
the composition series of any uniserial module.  On the other hand, condition
(N) is satisfied by a class of algebras within which {\it arbitrary} affine
varieties can still be realized as uniserial varieties $V_{\SS}$ (see
Proposition 3.8 for a precise statement).

In Section 4, we completely pin down the finite dimensional algebras for
which all of the varieties $V_{\SS}$ are finite.  This condition is stronger
than finite uniserial type since, in the presence of oriented cycles of
$\Gamma$, the canonical surjections $V_{\SS} \rightarrow \{$isomorphism types
of uniserials in $\la$-mod with composition series $\SS \}$ can have infinite
fibres;  in fact, in contrast to finite uniserial type, this condition is
not left-right-symmetric in general.  The algebras $\la$ for which all the
uniserial varieties are finite turn out to be precisely the ones for which
there is a one-to-one correspondence between the isomorphism types of
uniserial left $\la$-modules and their graphs, all of the latter being
edge paths in that case.  Moreover, these algebras are particularly easy
to recognize via quiver and relations;  namely, they are characterized by
the following strengthening of Condition (N):  Whenever $p : e \rightarrow
e'$ in $K\Gamma$ is a mast of positive length and $\alpha: e \rightarrow
e'$ an arrow, $p$ may differ from $\alpha$ only through an oriented cycle
based at the vertex $e'$.

In Section 5, which contains the deepest results of the paper, we let
the base field $K$ be infinite and assume Condition (N).  In the
presence of this condition, we explore the sequences $\SS = (S(0), \dots,
S(l))$ of simple left
$\la$-modules which do arise as composition series of uniserial left
$\la$-modules, but only of finitely many.  Our main strategy is to use 
properties of the variety
$V_{\SS}$ to pin down the possible patterns for sequences $\SS$, and to
analyze quiver and relations of $\la$ `at' the (necessarily unique) mast $p$
of the uniserial modules with composition series $\SS$.  Since a complete
description of the paths
$p$ going with sequences $\SS$ of simple modules as above leads to
combinatorial distinctions that obscure rather than clarify the emerging
picture, we content ourselves -- in our main theorems -- with describing
relatively simple conditions which are necessary or sufficient and which, in
fact, `almost' meet. We supplement these conditions with an algorithmic
procedure for recognizing algebras of finite uniserial type on the basis of
their quivers and relations, which decides the issue in all cases.  The most
readily quotable result already gives a rough outline of the picture to be
drawn:  Namely, if
$\SS$ is a sequence of simple modules as above and $d$ is the dimension of the
variety $V_{\SS}$, then the top simple module $S(0)$ recurs at least $d$
times among $S(1), \dots S(l)$.  In fact, the sequence $\SS$, or
equivalently, the corresponding mast
$p$, has to be made up of repeating segments (see Preview 5.3).  This
implication hinges on the fact that each
$K$-transcendental variable in the coordinate ring of $V_{\SS}$ is naturally
associated with a pair $(u,v)$ of initial segments of $p$ such that
$v$ starts with an oriented cycle and then `essentially' repeats $u$.  The
cleanest form of this repetition yields a sufficient condition
for the finiteness of the number of uniserial modules with composition
series $\SS$.  We refer to Preview 5.3 for an overview, since a
bit more terminology is needed to formulate the results with precision. 
Several examples are given to demonstrate that the invariants of $\SS$
involved in our theory are readily accessible if $\la$ is given by means of
quiver and relations, and that these invariants permit a tight grip on the
uniserial representation theory via rather simple computations.

The most important invariants of a sequence $\SS$ -- or equivalently of a
mast
$p \in K\Gamma$ -- along this line are the slack halyards which are
introduced and discussed in Section 2.  Sections 6 and 7, finally, apply the
preceding theorems towards convenient characterizations of finite uniserial
type within restricted classes of algebras.  Section 6 focuses on monomial
relation algebras, Section 7 on algebras of Loewy length at most $6$;  in
the latter case our theory allows us to completely list the uniserial
modules over algebras of finite uniserial type.   

\subhead Setting and Prerequisites\endsubhead Throughout, we assume that $\la$
is a split basic algebra over an {\it infinite} field $K$, i.e., $\la$ can be
`coordinatized' in the form $\la\cong K\Gamma/I$, where $\Gamma$ is a
quiver and $I$ an admissible ideal in the path algebra $K\Gamma$. For
simplicity, we suppose that $\la=K\Gamma/I$ and identify the vertices of
$\Gamma$ with a privileged set of primitive idempotents of $\la$. Our
convention for composing paths in $K\Gamma$ will be as follows: `$qp$' stands
for `$q$ after $p$'. A {\it right subpath} of a path $p\in K\Gamma$ is a
path $u$ such that $p=u'u$ in $K\Gamma$ for some path $u'$; {\it left
subpaths} are defined symmetrically. Both right and
left subpaths are special cases of what we call {\it subpaths of $p$},
namely, paths $u$ such that $p=u'uu''$ for suitable paths $u'$ and $u''$ in
$K\Gamma$.

Given a left $\la$-module $M$, we call any element $m\in M\setminus JM$ with
the property that $m=em$ for one of the vertices $e$ of $\Gamma$ a {\it top
element} of $M$. Here $J$ denotes the Jacobson radical of $\la$. Recall that
given a uniserial module $U$ together with a top element $x\in U$, a great
deal of structural information about the module $U$ can be communicated by
means of its `layered and labeled graph' relative to $x$ [6, Section 2].
We informally remind the reader of such graphs: That $U$ be represented by a
graph of the form
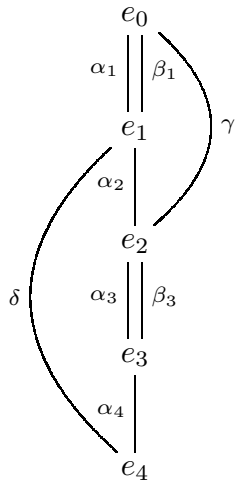
\ignore{
$$\xy\xymatrix{
e_0 \ar@{-}[d]<-0.5ex>_{\alpha_1} \ar@{-}[d]<0.5ex>^{\beta_1}
\ar@{-}@/^2pc/[dd]^\gamma\\
e_1 \ar@{-}[d]_{\alpha_2} \ar@{-}@/_2.75pc/[ddd]_\delta\\
e_2 \ar@{-}[d]<-0.5ex>_{\alpha_3} \ar@{-}[d]<0.5ex>^{\beta_3}\\
e_3 \ar@{-}[d]_{\alpha_4}\\
e_4}\endxy$$
}
\noindent relative to the top element $x$, means that the $\alpha_i$ are
arrows $e_{i-1}\rightarrow e_i$ for $1\le i\le4$, the $\beta_i$ are arrows
$e_{i-1}\rightarrow e_i$ for $1=1,3$, and $\gamma : e_0\rightarrow e_2$ and
$\delta : e_1\rightarrow e_4$ are arrows such that 
$$\gather \la \alpha_1x =\la\beta_1x= JU\\
\la\alpha_2\alpha_1x =\la\alpha_2\beta_1x =\la\gamma x =J^2U\\
\la\alpha_3\alpha_2\alpha_1x =\la\beta_3\alpha_2\alpha_1x =J^3U\\
\la\alpha_4\alpha_3\alpha_2\alpha_1x= \la\delta\alpha_1x =J^4U;\endgather$$
moreover, the graph tells us that the composition length of $U$ is 5 and
that, whenever $\epsilon : e_i\rightarrow e_j$ is an arrow in $\Gamma$
which does not correspond to a labeled edge $e_i\ {\vrule height2.5pt
depth-2.25pt width0.6truecm}\ e_j$ in the graph of
$U$, we have $\epsilon J^{i-1}x \subseteq J^{i+1}U$.

The varieties mentioned in the introduction are most readily manageable
when described in the framework of a fixed coordinatization of $\la$.
Recall that a {\it mast} of a uniserial module $U$ of length $l+1$ is any
path $p\in K\Gamma$ of length $l$ with the property that $pU\ne0$. (In the
previous example, the masts of $U$ are $\alpha_4\alpha_3\alpha_2\alpha_1$,
$\alpha_4\alpha_3\alpha_2\beta_1$, $\alpha_4\beta_3\alpha_2\alpha_1$,
$\alpha_4\beta_3\alpha_2\beta_1$.) Observe that, when the algebra $\la$ is
of finite uniserial type, its quiver $\Gamma$ is without double arrows,
which in turn implies that each uniserial $\la$-module has a unique mast.

Given a path $p\in K\Gamma$, we use the term {\it detour on $p$} for any
pair $(\alpha,u)$, where $u$ is a right subpath of $p$ (which may have
length zero) and $\alpha$ is an arrow such that $\alpha u$ is a path in
$K\Gamma$ which is not a right subpath of $p$ but has the following
property: there exists a right subpath $v$ of $p$, strictly longer than $u$,
which ends in the same vertex as $\alpha$ (see the figure below). Write
$(\alpha,u)\detour p$ to indicate that $(\alpha,u)$ is a detour on
$p$.
\ignore{
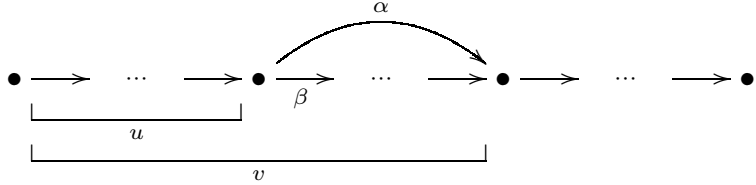
$$\xymatrixcolsep{1.5pc}
\xy\xymatrix{
\bullet \ar[r] \ar@^{|-|}[rrr]<-3ex>_u \ar@^{|-|}[0,6]<-6ex>_v
&\ar@{}[r]|\cdots &\ar[r]
&\bullet \ar[r]_\beta \ar@/^1.5pc/[rrr]^\alpha &\ar@{}[r]|\cdots
&\ar[r] &\bullet \ar[r] &\ar@{}[r]|\cdots &\ar[r] &\bullet
}\endxy$$
}
\noindent In case $(\alpha,u)$ is a detour on $p$, the family of
\underbar{all} right subpaths of $p$ which are longer than $u$ and end in
the same vertex as $\alpha$ is of key importance in the construction of
potential uniserial modules $U$ with mast $p$; this family is labeled
$(v_i(\alpha,u))_{i\in I(\alpha,u)}$. Indeed, for any such module $U$ --
with top element x, say -- the product $\alpha ux$ is a $K$-linear
combination of the elements $v_i(\alpha,u)x$, $i\in I(\alpha,u)$, and $U$
is determined up to isomorphism by $x$ and the coefficients arising in
these linear combinations. 

As shown in [6, Section 3], there is an affine algebraic variety $V_p$ (not
necessarily irreducible), the variables of which are in one-to-one
correspondence with the triples 
$(\alpha,u,i)$, where $i\in I(\alpha,u)$ and $(\alpha,u)\detour p$,
together with a natural surjection 
$$\Phi_p : V_p \rightarrow \{\text{isomorphism classes of uniserials in\ }
\lamod \text{\ with mast\ } p \}$$
which pins down the correlations among the parameters describing the
uniserial modules with mast $p$. This variety $V_p$ is nonempty if and only
if $p$ does arise as a mast of a uniserial module. Its coordinate ring will
be denoted by $A_p$.

In a few arguments of the present paper, we will refer to the construction
of a set of polynomials for $V_p$ in the commutative ring $K[X_i(\alpha,u)
\mid i\in I(\alpha,u),\ (\alpha,u)\detour p]$. Such polynomials can be
obtained from a representative set of relations in $K\Gamma$ by reducing
the relations in this set modulo a certain equivalence relation `$\seq$' on
the enlarged polynomial ring $K\Gamma[X_i(\alpha,u)
\mid i\in I(\alpha,u),\ (\alpha,u)\detour p]$ which is defined by what we
call the `substitution equations for $p$'. The reader can find full detail
on this point, as well as on the basic properties of $V_p$ and $\Phi_p$, in
[6, Section 3].

In [6], it is also explained how to move from the coordinatized
varieties $V_p$ to a family $V_\SS$ of varieties which describes all the
uniserial modules having a fixed sequence $\SS= (S(0),\dots,S(l))$ of
consecutive composition factors, up to isomorphism (by a slight abuse of
language, we will also speak of uniserial modules with `composition series
$\SS$'). In case $\Gamma$ is
without double arrows -- this being the case which interests us primarily,
since it subsumes finite uniserial type -- this transition is particularly
straightforward. Namely, if
$(e(0),\dots,e(l))$ is the sequence of vertices corresponding to the simple
modules $S(i)$ in $\SS$, we then have at most one path $p$ of length $l$
passing through these vertices in the given order, and the image of the map
$\Phi_p$ is the set of all isomorphism classes of uniserial left
$\la$-modules with composition series
$\SS$. Consequently, writing $V_\SS$ instead of $V_p$ and
$\Phi_\SS$ instead of
$\Phi_p$ will not cause any ambiguities; to see that this change
of notation is meaningful, recall that the isomorphism type of the variety
$V_\SS$ is not affected by the coordinatization of $\la$ (see [7]).

For the basics of algebraic geometry and commutative algebra, we refer the
reader to the texts by Fulton [3], Hartshorne [5], and Matsumura [8].

\head 2. Halyards -- the role of circular halyards\endhead

The purpose of this section is to introduce `halyards' on a mast $p$, and to
discuss some of their basic properties. Roughly speaking, halyards are detours
on
$p$ which make an appearance in at least one graph of a uniserial module with
mast $p$. They will turn out to be crucial in the process of
determining algebras of finite uniserial type. All of the following concepts
are, a priori, tied to the chosen coordinatization of
$\la= K\Gamma/I$. However, when $\Gamma$ is without double arrows -- in
particular, when $\la$ is of finite uniserial type -- this problem
disappears.

\definition{Definitions 2.1} Suppose that $p\in K\Gamma$ is a mast of a
uniserial left $\la$-module, and let $(\alpha,u)$ be a detour on $p$.
Moreover, let
$(v_i(\alpha,u))_{i\in I(\alpha,u)}$ be the family of right subpaths of $p$
longer than $u$ which end in the same vertex as $\alpha$, and let
$X_i(\alpha,u)$ be the corresponding variables. As before, $A_p$ stands for
the coordinate ring of the variety $V_p$. 

(1) The detour $(\alpha,u)$ is called a {\it halyard on $p$} if at least one
of the variables $X_j(\alpha,u)$ has nonzero image in $A_p$, that is, if there
exists a uniserial
$\la$-module
$U$ with mast
$p$ and top element
$x$ such that
$\alpha ux\ne 0$. Equivalently, this means that the graph of some uniserial
module with mast $p$ has a subgraph of the form
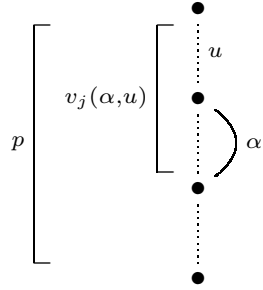
\ignore{
$$\xymatrixrowsep{1.5pc}
\xy\xymatrix{
\bullet \ar@{.}[d]^u \ar@^{|-|}[dd]<-3ex>_{v_j(\alpha,u)}
\ar@^{|-|}[ddd]<-12ex>_p\\
\bullet \ar@{.}[d] \ar@{-}@/^1pc/[d]^\alpha\\
\bullet \ar@{.}[d]\\
\bullet}\endxy$$
}
\noindent In this situation, we also say that the uniserial module $U$ with
mast $p$ has halyard $(\alpha,u)$, or that the above is a graph of the
halyard $(\alpha,u)$ on $p$.

(2) Suppose that $(\alpha,u)$ is a halyard on $p$ and $j\in I(\alpha,u)$.
Then $X_j(\alpha,u)$ is called a {\it slack} variable of $(\alpha,u)$ in
case the image of $X_j(\alpha,u)$ in $A_p$ is transcendental over $K$. A
halyard with at least one slack variable is in turn called {\it slack}.
Non-slack halyards (and variables) will be called {\it tight}.

(3) A halyard $(\alpha,u)$ on a mast $p$ is called {\it circular} if it has
a graph of the form
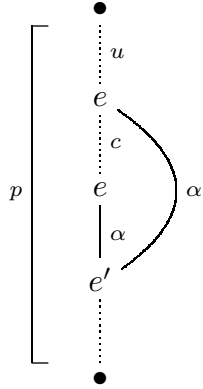
\ignore{
$$\xymatrixrowsep{1.5pc}
\xy\xymatrix{
\bullet \ar@{.}[d]^u \ar@^{|-|}[dddd]<-5ex>_p\\
e \ar@{.}[d]^c \ar@{-}@/^2pc/[dd]^\alpha\\
e \ar@{-}[d]^\alpha\\
e' \ar@{.}[d]\\
\bullet}\endxy$$
}
\noindent where $c : e\rightarrow e$ is a non-trivial oriented
cycle.\enddefinition

According to the following lemma, a mast that has only tight halyards will
support only finitely many uniserial modules.

\proclaim{Lemma 2.2} Let $(\alpha,u)$ be a halyard on a mast $p$. A variable
$X_j(\alpha,u)$ is slack if and only if the canonical projection
$V_p\rightarrow K$, $(k_i(\beta,v)) \mapsto k_j(\alpha,u)$, has infinite
image. The latter is in turn equivalent to the following
representation-theoretic condition: There are infinitely many elements $k\in
K$ with the property that there exists a uniserial module $U_k\in\lamod$
having mast $p$ and a top element $x_k$ with $\alpha ux_k
\equiv kv_j(\alpha,u)x_k$ $(\text{mod\ } \sum_{i\in I(\alpha,u),\ i\ne j}
Kv_i(\alpha,u)x_k)$.
\qed\endproclaim

Loosely, the reason for our interest in circular halyards is as
follows: As we wish to identify those slack halyards which entail the
existence of infinitely many uniserial modules, we are in particular
interested in those halyards which are responsible for failure of
bijectivity of the maps $\Phi_p$ from the varieties $V_p$ to the
corresponding sets of isomorphism types of uniserials with mast $p$. The
culprits are slack circular halyards.

\example{Remarks and Examples 2.3} (a) As an immediate consequence of Lemma
2.2, we obtain the following implication: In case there exist infinitely
many pairwise non-isomorphic uniserial left
$\la$-modules with mast $p$, there is a slack halyard on $p$. However, even
over an infinite field $K$, the converse fails in general. Indeed,
if
$\la= K\Gamma/
\langle
\beta^2\rangle$, where $\Gamma$ is the quiver
\ignore{
$$\xymatrixcolsep{3pc}
\xy\xymatrix{
1 \lloopd{\beta} \ar[r]^\alpha &2
}\endxy$$
}
\noindent then it is readily checked that $p=\alpha\beta$ is a mast with
slack halyard $(\alpha,e)$; on the other hand, there is precisely one
uniserial left $\la$-module with mast $p$, up to isomorphism (see [6,
Section 4]).

In particular, there will be only finitely many isomorphism classes of
uniserials with mast $p$ in $\lamod$, whenever we can show that the
slack halyards on $p$ are irrelevant to the isomorphism types of the
corresponding uniserial modules in the following sense: Suppose that
$(\alpha_1,u_1),\dots,(\alpha_m,u_m)$ are the slack halyards on $p$ and
$I_s(\alpha_i,u_i)\subseteq I(\alpha_i,u_i)$ the set of indices of the slack
variables $X_j(\alpha_i,u_i)$. Then the
number of uniserial left $\la$-modules with mast $p$ is finite provided that
each such uniserial module $U$ has a top element $x$ with the property that
$$\alpha_iu_ix\ \in\ \sum_{j\in I(\alpha_i,u_i)\setminus I_s(\alpha_i,u_i)}
Kv_j(\alpha_i,u_i)x$$
for all $i=1,\dots,m$.

 (b) If $(\alpha,u)$ is a circular halyard on a
mast $p$ and $c$ an oriented cycle as in the definition, then $c$
necessarily starts in an arrow $\beta\ne \alpha$ because, in particular,
$(\alpha,u)$ is a detour on $p$. A more detailed rendering of a graph thus
looks as follows:
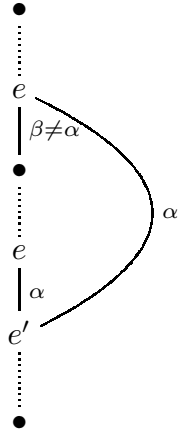
\ignore{
$$\xymatrixrowsep{1.25pc}
\xy\xymatrix{
\bullet \ar@{.}[d]\\
e \ar@{-}[d]^{\beta\ne\alpha} \ar@{-}@/^3.5pc/[ddd]^\alpha\\
\bullet \ar@{.}[d]\\
e \ar@{-}[d]^\alpha\\
e' \ar@{.}[d]\\
\bullet}\endxy$$
}

Obviously the circular halyards $(\alpha,u)$ on $p$ are exactly the ones
which disappear on dualizing -- more precisely, they are the halyards on $p$
which do not show up in any graph of a right $\la$-module $D(U)$, where $U$
runs through the uniserial left $\la$-modules with mast $p$; here $D$
denotes the standard duality from $\lamod$ to mod-$\la$.

(c) Suppose that $p= \alpha c : e\rightarrow e'$ is a mast, where $\alpha :
e\rightarrow e'$ is an arrow and  $c : e\rightarrow e$  an oriented
cycle of positive length which does not pass through the vertex $e'$; in
particular,
$e\ne e'$. Then
$(\alpha,e)$ is a circular halyard on $p$ having graph
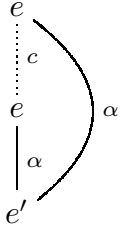
\ignore{
$$\xymatrixrowsep{1.8pc}
\xy\xymatrix{
e \ar@{.}[d]^c \ar@{-}@/^2pc/[dd]^\alpha\\
e \ar@{-}[d]^\alpha\\
e'}\endxy$$
}
\noindent By Lemma 3.1, this halyard is
necessarily slack.

Note that $(\alpha,e)$ need not be a circular halyard in case the vertex $e'$
lies on the cycle
$c$, for then it can occur that all graphs of uniserials with mast $p$ are
of the form
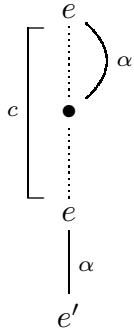
\ignore{
$$\xymatrixrowsep{1.8pc}
\xy\xymatrix{
e \ar@{.}[d] \ar@{-}@/^1pc/[d]^\alpha \ar@^{|-|}[dd]<-3ex>_c\\
\bullet \ar@{.}[d]\\
e \ar@{-}[d]^\alpha\\
e'}\endxy$$
}
\noindent For instance, consider the algebra
$\la= K\Gamma/I$ with quiver $\Gamma$
\ignore{
$$\xymatrixcolsep{3pc}
\xy\xymatrix{
1 \lloopd{\alpha} \rloopd{\beta}
}\endxy$$
}
\noindent and $I=\langle \beta\alpha -\beta^2,\ \alpha\beta -\alpha^2,\
\text{all paths of length 3}\rangle$. While $p=\alpha\beta$ is a
mast and $\beta=c$ an oriented cycle, there are clearly no circular halyards
on
$p$. In fact, there are no circular halyards on {\it any} mast of a uniserial
left
$\la$-module in this example.

(d) The following is an example of a slack halyard with two variables, one
of which is tight. This time, let $\la= K\Gamma/ \langle \gamma^2,\
\gamma\beta_1\alpha- \gamma\beta_2\alpha \rangle$, where $K$ is an infinite
field and
$\Gamma$ is the quiver
\ignore{
$$\xymatrixcolsep{3pc}
\xy\xymatrix{
1 \ar[r]^\alpha &2 \ar[r]<0.5ex>^{\beta_1}
\ar[r]<-0.5ex>_{\beta_2} &3 \rloopd{\gamma}
}\endxy$$
}
\noindent Then $p=\gamma\beta_1\alpha$ is a mast with a
single halyard, namely $(\beta_2,\alpha)$, whose variables $X_1=
X_1(\beta_2,\alpha)$ and $X_2= X_2(\beta_2,\alpha)$ correspond to the right
subpaths $v_1= \beta_1\alpha$ and $v_2=p$ of $p$. By [6, Section 3], $V_p=
V(X_1-1)\subseteq \AA^2$, i.e., $X_1$ is a tight variable of the halyard
$(\beta_2,\alpha)$, whereas $X_2$ is slack.

(e) Finally, we exhibit tight halyards with variables assuming more than a
single value on $V_p$. Consider the algebra $\la= K\Gamma/I$, where
$\text{char} K\ne 2$, the quiver
$\Gamma$ is
\ignore{
$$\xymatrixcolsep{3pc}
\xy\xymatrix{
1 \ar[r]<0.5ex>^{\alpha_1} \ar[r]<-0.5ex>_{\beta_1} &2
\ar[r]<0.5ex>^{\alpha_2}
\ar[r]<-0.5ex>_{\beta_2} &3
}\endxy$$
}
\noindent and $I=\langle \alpha_2\alpha_1 -\beta_2\beta_1,\ \beta_2\alpha_1
-\alpha_2\beta_1 \rangle$. If $p= \alpha_2\alpha_1$, then $V_p= V(X_1X_2-1,\
X_2-X_1) \subseteq \AA^2$, i.e., $V_p= \{(1,1),(-1,-1)\}$. In particular,
 the detours $(\beta_1,e_1)$ and $(\beta_2,\alpha_1)$ on $p$ are both tight
halyards on $p$. \qed\endexample

\subhead Coordinate-free description of slack halyards\endsubhead 

Let us focus on the uniserial left $\la$-modules with a fixed sequence $\SS
= (S(0),\dots,S(l))$ of consecutive simple composition factors. As is readily
seen, the existence of {\it tight} halyards of such uniserial modules may
depend on the choice of coordinate system. On the other hand, the existence
of {\it slack} halyards of such modules is independent of that choice.

We specialize to the situation where $\Gamma$ is without double arrows, which
presents us with a slightly simplified picture. Let $(e(0),\dots,e(l))$ be
the sequence of vertices corresponding to the simple modules $S(i)$ in
$\SS$, and let $p$ be the unique path of length $l$ passing through these
vertices in the given order.

Suppose that
$p$ is a mast or, equivalently, that $V_\SS$ is nonempty. A halyard
$(\alpha,u)$ on $p$ can be communicated by a pair $(S,\rho)$, where
$\rho=\len(u)$ indicates the radical layer from which the halyard emerges
in any of its graphs, and $S$ represents the isomorphism class of the simple
module $\la e/Je$ corresponding to the terminal vertex $e$ of $\alpha$. By
the definition of a detour, $S$ is then isomorphic to at least one of the
$S(i)$ with
$\rho+1\le i\le l$. Clearly, the right subpaths
$v_i(\alpha,u)$,
$i\in I(\alpha,u)$, of $p$ are in 1--1 correspondence with those indices
$\sigma\in \{\rho+1, \dots,l\}$ for which $S(\sigma)\cong S$; in other
words, the variables $X_i(S,\rho)= X_i(\alpha,u)$ of the halyard
$(S,\rho)$ may be identified with the lengths of the corresponding right
subpaths
$v_i(\alpha,u)$. The reader will readily verify that the slack variables
$\sigma$ of the pair $(S,\rho)$ are independent of the coordinatization. It
is thus unequivocal to speak of a `slack variable
$\sigma$ of a halyard $(S,\rho)$ on a sequence $\SS$' of simple left
$\la$-modules. For instance,
in the example preceding these remarks, $V_\SS= \AA^1$ for $\SS=
(S(0),S(1),S(2),S(3))= (S_1,S_2,S_2,S_3)$, and the pair $(S_3,1)$ is a slack
halyard on $\SS$, with slack variable 3.

Finally, note the following immediate consequence of our definitions: There
are no slack halyards on the sequence
$\SS$ of simples if and only if $V_\SS$ is
finite.

\head 3. A strong necessary condition for finite uniserial type\endhead

We start our investigation of finite uniserial type by establishing a
necessary condition (N). This condition narrows down the possibilities
for the quiver
$\Gamma$ and the relations inasmuch as it imposes severe restrictions on those
paths which arise as masts. In particular, (N) implies that all halyards,
slack or tight, on masts of uniserial modules are circular. Consequently,
this condition is tantamount with finite uniserial type of $\la$ whenever
the quiver 
$\Gamma$ is acyclic.  We will derive further consequences from
(N), in order to demonstrate to what degree this condition impinges on the
isomorphism classes and graphs of the uniserial representations in general.

On the other hand, we will indicate how complex the uniserial representation
theory of algebras satisfying (N) can still be, by showing that there is
enough leeway to realize any affine algebraic variety as a variety $V_p$
over such an algebra, even if we impose the additional requirement that the
corresponding canonical map from $V_p$ to the set of isomorphism classes of
uniserials with mast $p$ be bijective.

\proclaim{Lemma 3.1} Let $\alpha : e\rightarrow e'$ be an arrow, where $e$ and
$e'$ are (not necessarily distinct) vertices of $\Gamma$. Moreover, suppose
that $p : e\rightarrow e'$ is a mast of positive length which does not start
in $\alpha$, i.e., $p=p'\beta$ for some arrow $\beta\ne \alpha$. Then
$(\alpha,e)$ is a slack halyard on $p$.

More precisely: If $v_1,\dots,v_s=p$ are the subpaths of positive length of
$p$ which end in $e'$, then for each $k\in K$, there exists a
uniserial module $U_k\in \lamod$ with top element $z_k$ such that 
$$\alpha z_k= k_1v_1z_k +\dots+ k_{s-1}v_{s-1}z_k +kpz_k$$
for suitable $k_i\in K$; in particular, the variable $X_s(\alpha,e)$ going
with $v_s=p$ is slack.\endproclaim

\demo{Proof} Recall that we obtain $V_p$ by inserting the substitution
equation
$$\alpha \seq X_1v_1 +\dots+ X_sv_s = X_1(\alpha,e)v_1 +\dots+
X_s(\alpha,e)v_s$$ 
from the right into a selection of relations from $I$ (see
[Geom.I, Section 3]). Say $q\alpha$ is a path occurring non-trivially in such
a relation. Then
$\len(q)\ge 1$ and, in the substitution process, $q\alpha$ is replaced by
the element $X_1qv_1 +\dots+ X_sqv_s \in K\Gamma[X]$. Since $v_s=p$ and
therefore $\len(qv_s) >\len(p)$, the path $qv_s$ is not a route on $p$, and
consequently, the term $X_sqv_s$ is set equal to zero in the next step of
the substitution process. Thus the variable $X_s$ drops out of all of the
polynomials in $K[X]$ resulting from the substitution process, irrespective of
the relations in $I$ with which we start this process. Since, on one hand,
$V_p$ is the vanishing locus of a set of polynomials obtained from elements
of $I$ via iterated substitution, and since, on the other hand
$V_p\ne \varnothing$, the projection of
$V_p$ onto the coordinate $X_s= X_s(\alpha,e)$ equals $K$. \qed\enddemo

\proclaim{Proposition 3.2} Let $\alpha : e\rightarrow e'$ be an arrow and $p :
e\rightarrow e'$ in $K\Gamma$ a mast of positive length which neither starts
in $\alpha$ nor ends in $\alpha$. Then there exists a subpath $q$ of $p$ with
the property that there are infinitely many isomorphism types of uniserial
left $\la$-modules with mast $q$.\endproclaim

\demo{Proof} We will prove the proposition by
induction on $\len(p)$. Again let $v_1,\dots,v_s=p$ be all the right subpaths of
positive length of $p$ which end in $e'$. 

If $\len(p)=1$, that is, if $p$ is an arrow from $e$ to $e'$, the fact that
$\alpha\ne p$ clearly yields infinitely many pairwise nonisomorphic
uniserials with graph
\ignore{
$$\xymatrixrowsep{1.5pc}
\xy\xymatrix{
e \ar@{-}[d]_p \ar@{-}@/^1pc/[d]^\alpha\\
e'
}\endxy$$
}
\noindent (just keep in mind that $K$ is infinite).

Now suppose that $\len(p)\ge 2$, say $p=p'\beta$ for some arrow $\beta\ne
\alpha$, and that for any arrow $\gamma : e_0\rightarrow e_0'$ and any mast
$q : e_0\rightarrow e_0'$ with $0<\len(q)<\len(p)$ our claim is true.
Moreover, assume the proposition to be false for the mast $p$.

From Lemma 3.1, we know that $(\alpha,e)$ is a halyard on $p$ and that the
variable $X_s=X_s(\alpha,e)$ in the corresponding substitution equation
$$\alpha \seq X_1v_1 +\dots+ X_sv_s$$
is slack. Our assumption thus gives us distinct scalars $k=k_s,\ l=l_s\in
K$, together with isomorphic uniserial modules $U_k$ and $U_l$ having mast
$p$ and top elements $x$ and $y$, respectively, such that
$$\align \alpha x &= \sum_{i=1}^{s-1} k_iv_ix + kpx\\
\alpha y &= \sum_{i=1}^{s-1} l_iv_iy + lpy\endalign$$
for suitable scalars $k_i,l_i$, $i\le s-1$. Let $w_1,\dots,w_t$ be all the
right subpaths of $p$ of positive length which end in the vertex $e$, listed
in order of increasing length. It is clearly harmless to assume that
$U_k=U_l=U$ and $y=x+ \sum_{j=1}^t c_jw_jx$ for certain $c_j\in K$. Expanding
both sides of the second equation by inserting the first, we thus obtain
$$\sum_{i=1}^s (k_i-l_i)v_ix = \sum_{i,j}
l_ic_jv_iw_jx - \sum_{j=1}^t c_j\alpha w_jx.$$
Since the $v_ix$ are $K$-linearly independent in $U$, the nonzero term
$(k-l)px$ on the left-hand side does not cancel against other terms on the
left.

 Let us next inspect the terms $\alpha w_jx$ occurring on the
right-hand side: For $1\le j\le t$, we let $w_j' : e\rightarrow e'$ be the
left subpath of $p$ such that $p= w_j'w_j$, this last equality holding in
$K\Gamma$. In view of the fact that $p$ does not end in $\alpha$, none of
the left subpaths $w_j'$ ends in $\alpha$. By assumption, there are only
finitely many uniserial left $\la$-modules, up to isomorphism, having masts
among the subpaths of
$w_1',\dots,w_t'$, and therefore our induction hypothesis forces those of
the $w_j'$ which have positive length to start in $\alpha$. The only $w_j'$
which is potentially of length zero is $w_t'$; in case it is, we conclude
$e=e'$ and $w_t=p$, which gives us $\alpha w_tx=0$. Thus, whenever $\alpha
w_jx\ne 0$, the path $\alpha w_j$ is a right subpath of $p$ in $K\Gamma$.
None of these paths $\alpha w_j$ is equal to $p$ since, by hypothesis, $p$
does not end in $\alpha$. Again using the fact that the elements $qx$, where
$q$ runs through the right subpaths of $p$, are $K$-linearly independent, we
infer that $v_iw_jx\ne 0$ for some $i$ and $j$. Since this inequality
forces the length of the path $v_iw_j$ to be bounded above by that of $p$,
we derive $i\le s-1$ and $\len(w_j)< \len(p)$; in particular, the above
considerations show that $\alpha w_j$ is a right subpath of $p$. Being a
right subpath of positive length of $p$, the path $v_i$ starts in the same
arrow $\beta\ne \alpha$ as $p$, and from $v_iw_jx\ne 0$ we deduce $\beta
w_jx\ne 0$. In other words, the detour $(\beta,w_j)$ on $p$ is actually a
halyard showing up in the graph of $U$ relative to $x$ as follows:
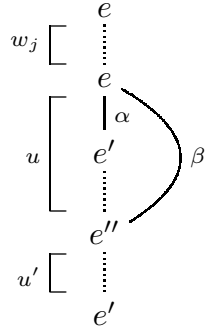
\ignore{
$$\xymatrixrowsep{1pc}
\xy\xymatrix{
e \ar@{.}[d] \ar@^{|-|}[d]<-4ex>_{w_j}\\
e \ar@{-}[d]^\alpha \ar@{-}@/^2pc/[dd]^\beta
\ar@^{|-|}[dd]<-4ex>_u\\ 
e' \ar@{.}[d]\\
e'' \ar@{.}[d] \ar@^{|-|}[d]<-4ex>_{u'}\\
e'
}\endxy$$
}
\noindent Here $e''$ is the terminal vertex of $\beta$, $u$ is the subpath
of $p$ with the property that $\la \beta w_jx =\la uw_jx$, and $u'$ is the
left subpath of $p$ such that $p= u'uw_j$. We apply the induction
hypothesis again, this time to the mast $u : e\rightarrow e''$ and the arrow
$\beta : e\rightarrow e''$. Since by assumption there are only finitely many
isomorphism types of uniserial modules having as mast a subpath of $u$, and
since $u$ does not start in $\beta$, the last arrow of $u$ must be $\beta$.
This means that
$u'\beta :e\rightarrow e'$ is a left subpath of $p$ and, in particular,
again a mast. Clearly, 
$$0< \len(u'\beta) < \len(p)-\len(w_j) <\len(p),$$
so, combining once more our assumption with the induction hypothesis --
applied to the mast $u'\beta :e\rightarrow e'$ and the arrow $\alpha$ -- we
infer that $u'\beta$ ends in $\alpha$. But $u'\beta$ being a left subpath
of $p$, this implies that $p$ ends in $\alpha$, a contradiction to the
hypothesis. \qed\enddemo

\proclaim{Theorem 3.3} {\rm (I)} If $\la$ has finite uniserial type, then the
following condition (N) is satisfied: Whenever
$\alpha : e\rightarrow e'$ is an arrow in $\Gamma$ and $p : e\rightarrow e'$
a mast of positive length, then $p\in K\Gamma\alpha \cup \alpha K\Gamma$, that
is,
$p$ is of the form
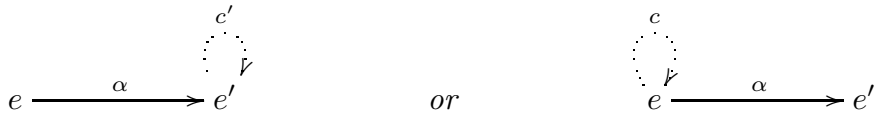
\ignore{
$$\xy\xymatrix{
e \ar[rr]^\alpha &&e' \udotloopr{c'} &&\text{or}&& e \ar[rr]^\alpha
\udotloopr{c} &&e'
}\endxy$$
}
\noindent where $c',c$ are oriented cycles which may be trivial.

{\rm (II)} Condition (N), in turn, implies that:

{\rm (a)} $\Gamma$ has no double arrows. Hence, each uniserial left
$\la$-module has a unique mast, and all the varieties $V_\SS$, where $\SS$
is a sequence of simple $\la$-modules, are uniquely determined by the
isomorphism type of $\la$.

{\rm (b)} All halyards on masts of uniserial left $\la$-modules are
circular. In particular, given any repetition-free sequence $\SS$ of simple
objects in $\lamod$, there is at most one uniserial module with composition
series $\SS$.

{\rm (c)} If $\gamma$ is a loop attached to a vertex $e$ of $\Gamma$ and $p :
e\rightarrow e$ a mast, then $p= \gamma^m$ for some $m$, or, in
coordinate-free terms: If there exists
a uniserial left $\la$-module $W$ of length 2 with $\top W\cong \soc W\cong
S$ (i.e., if $\Ext(S,S)\ne 0$), and if $U$ is any uniserial left $\la$-module
of length
$l\ge2$ with $\top U\cong \soc U\cong S$, then $J^iU/J^{i+2}U\cong W$ for
all $i\le l-2$.

{\rm (d)} Suppose that $\SS= (S(0),S(1),\dots,S(l))$ is a sequence of
simple left $\la$-modules such that $V_{(S(1),\dots,S(l))}$ is finite. Then
there are only finitely many uniserial left
$\la$-modules with composition series $\SS$.\endproclaim

\remark{Remark} A coordinate-free rendering of condition (N) is as follows:
If there exists a uniserial left $\la$-module $W$ of length 2 with top $S$
and socle $S'$, and if $U$ is any uniserial left $\la$-module of length
$l\ge2$ with top $S$ and socle
$S'$, then either $U/J^2U\cong W$ or else $J^{l-2}U\cong W$.\endremark

\demo{Proof} Part (I) follows immediately from Proposition 3.2, and all
but the last of the statements under (II) are clear.

So suppose that (N) holds, and that $V_\SS$ is nonempty. Let $p :
e(0)\rightarrow e(l)$ in $K\Gamma$ be
the unique path of length $l$ passing through the vertices that correspond
to the simple modules $S(i)$ in the given order. Because of our
hypothesis on $V_{(S(1),\dots,S(l))}$, we need only show that halyards of the
form $(\beta,e(0))$ on
$p$ do not affect the isomorphism type of the corresponding uniserial module
(cf\. Remark 2.3(a)). Suppose that $(\beta_1,e(0)),\dots, (\beta_m,e(0))$ is a
complete list of the distinct halyards of this type on the mast $p$, and let
$t(i)$ be the terminal vertex of
$\beta_i$. Since $\Gamma$ has no double arrows, the list
$t(1),\dots,t(m)$ is without repetitions. Moreover, for each $i\in
\{1,\dots,m\}$, we list the different right subpaths of
$p$ of positive length ending in $t(i)$, say $v_{i1},\dots,v_{is_i}$. Since
$p$ does not start in any of the $\beta_i$, condition (N) tells us that each
$v_{ij}$ ends in
$\beta_i$, say $v_{ij} =\beta_iw_{ij}$, where $w_{ij} : e(0)\rightarrow e(0)$
is an oriented cycle of positive length.

Now let $U$ be a uniserial module with mast $p$ and top element $x$ such
that $\beta_ix =\sum_{j=1}^{s_i} k_{ij}v_{ij}x$ for $1\le i\le m$. We only
need to find a top element $y$ of $U$ such that $\beta_iy=0$ for $1\le i\le
m$. 
We leave it to the reader to establish the existence of scalars $c_{ij}$ such
that the element $y= x+
\sum_{i\le m} \sum_{j\le s_i} c_{ij}w_{ij}x$ is as required. (The argument
is an easy special case of that given to prove Theorem 5.7.)
\qed\enddemo

Observe that Theorem 3.3 narrows down the sequences of simple
modules that potentially occur as sequences of consecutive composition factors
of uniserial left $\la$-modules in the case of finite uniserial type. For
example, if $S$ is a simple module with $\Ext(S,S)\ne 0$, then
composition factors isomorphic to $S$ can only occur in a single `connected
segment' within a composition series of a uniserial module $U$. The last
statement under (II) essentially says that, in the presence of condition (N),
halyards emerging from the top vertex of a mast do not affect the isomorphism
types of the pertinent uniserial modules. Moreover, since condition (N)
forces all halyards over an algebra of finite uniserial type to be circular,
we immediately obtain

\proclaim{Corollary 3.4} If $\Gamma$ is acyclic, then the following are
equivalent:

{\rm (1)} $\la$ has finite uniserial type.

{\rm (2)} $\la$ satisfies condition (N).

{\rm (3)} Whenever $S$ and $S'$ are simple left $\la$-modules
with $\Ext(S,S')\ne 0$, there exists precisely one uniserial
left $\la$-module $U$ with $\soc U\cong S'$ and $\top U\cong
S$. (In particular, there are no uniserial left $\la$-modules
of length $>2$ having socle $S'$ and top $S$.) \qed\endproclaim

Our next application of Theorem 3.3 concerns finite dimensional algebras
associated with tiled classical orders over discrete valuation domains. The
author was alerted to the homological interest of these algebras by E.
Kirkman and J. Kuzmanovich, and the following corollary is essentially due
to them. Let $F$ be a field with a discrete valuation, $D$ its valuation
domain, and $\pi$ a uniformizing parameter. A classical $D$-order $\O\subseteq
M_n(F)$ is called {\it tiled} if $\O$ contains a complete set of $n$
orthogonal idempotents. If we denote by $\la$ the finite dimensional algebra
$\O/\pi\O$ over the residue class field $K=D/(\pi)$ of $D$, then, by [4],
$$\lfd\O =1+\lfd\la,$$
where $\lfd\O$ denotes the supremum of the finite projective dimensions
attained on finitely generated left $\O$-modules. Moreover, it was observed
by Kirkman and Kuzmanovich that $\la=K\Gamma/I$, where $\Gamma$ is a quiver
obtained from the valued quiver  $\Gamma_{\text{val}}$ of $\O$ (see
[10]) via potential deletion of some loops, and
$I\subseteq K\Gamma$ is an admissible ideal with generators as follows: If
$p=\alpha_l\cdots\alpha_1$ is a path of length $l$ in $K\Gamma$, let $v(p)$
be the sum of the values of the arrows $\alpha_i$ in $\Gamma_{\text{val}}$.
Then $I$ is generated by the paths $p : e\rightarrow e'$ in $K\Gamma$ with
the property that $v(p)>v(q)$ for some path $q : e\rightarrow e'$, together
with the differences $p-q$ of paths sharing initial and terminal vertices such
that
$v(p)=v(q)$. Since, moreover, $\Gamma$ contains at least one path
$e\rightarrow e'$ for any two vertices $e$ and $e'$, it follows that each
simple left $\la$-module occurs with multiplicity precisely 1 as a composition
factor of each of the indecomposable projective modules $\la e$. It is easy to
derive (N) from these properties, as well as the fact that no uniserial left
$\la$-module has multiple composition factors. By Theorem 3.3(II),
this yields

\proclaim{Corollary 3.5} If $\O$ is a tiled classical order over a discrete
valuation ring with uniformizing parameter $\pi$, the algebra $\la=
\O/\pi\O$ has finite uniserial type.\qed\endproclaim

A simple example of an algebra which satisfies condition (N) but has
infinite uniserial type is as follows.

\example{Example 3.6} If $\Gamma$ is the quiver
\ignore{
$$\xymatrixcolsep{3pc}
\xy\xymatrix{
3 \ar[r]<-0.75ex>_\delta &1 \ar[l]<-0.75ex>_\gamma
\ar[r]<0.75ex>^\alpha &2 \ar[l]<0.75ex>^\beta
}\endxy$$
}
\noindent then $\la =K\Gamma/\langle \text{all paths of length\ }\ge
5\rangle$ satisfies (N). There is precisely one uniserial module with mast
$p=\delta\gamma\beta\alpha$; its graphs are
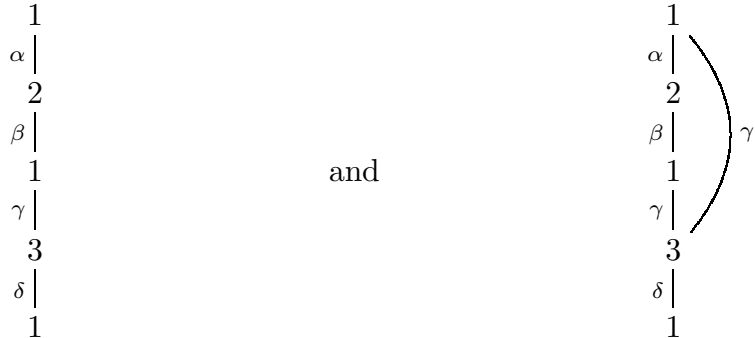
\ignore{
$$\xymatrixrowsep{1pc}
\xy\xymatrix{
1 \ar@{-}[d]_\alpha &&& &&& 1 \ar@{-}[d]_\alpha
\ar@{-}@/^1.5pc/[ddd]^\gamma\\
2 \ar@{-}[d]_\beta &&& &&& 2 \ar@{-}[d]_\beta\\
1 \ar@{-}[d]_\gamma &&& \text{and} &&& 1 \ar@{-}[d]_\gamma\\
3 \ar@{-}[d]_\delta &&& &&& 3 \ar@{-}[d]_\delta\\
1 &&& &&& 1
}\endxy$$
}
\noindent according to the choice of top element. However, there are
infinitely many uniserial modules with mast
$q=
\alpha\delta\gamma\beta$ having graph
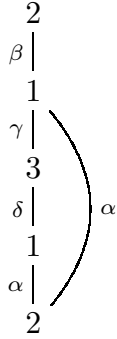
\ignore{
$$\xymatrixrowsep{1pc}
\xy\xymatrix{
2 \ar@{-}[d]_\beta\\
1 \ar@{-}[d]_\gamma \ar@{-}@/^1.5pc/[ddd]^\alpha\\
3 \ar@{-}[d]_\delta\\
1 \ar@{-}[d]_\alpha\\
2
}\endxy$$
}
\noindent In particular, $\la$ does not have finite uniserial type.
\qed\endexample

Without any further restrictions on $\Gamma$ and $I$, the uniserial
representation theory of an algebra $\la$ satisfying (N) may still be
arbitrarily involved. This statement is rendered more precise by the
following

\proclaim{Proposition 3.7} Given any affine algebraic variety $V$, there
exists a path algebra modulo relations, $\la= K\Gamma/I$, satisfying condition
(N), together with a sequence $\SS$ of simple left $\la$-modules such that
$V_\SS\cong V$. Moreover,
$\SS$ can be chosen so that the canonical map
$$\Phi_\SS : V_\SS \longrightarrow \{\text{isomorphism types of uniserials
in\ } \lamod\ \text{with composition series\ } \SS\}$$
is bijective.\endproclaim

\demo{Proof} As pointed out in the proof of Theorem G of [6], we do not
lose generality in assuming that $V=V(f_1,\dots,f_M)$ is the vanishing set
of polynomials $f_i\in K[X_1,\dots,X_m]$ with the property that each
variable $X_i$ occurs to a power $\le 1$ in any monomial involved in the
$f_i$. If $\P$ is the power set of $\{1,\dots,m\}$, we can thus write the
$f_i$ in the form $f_i= \sum_{A\in\P} c_i(A) \prod_{j\in A} X_j$ with
suitable scalars $c_i(A)\in K$.

Consider the quiver
\ignore{
$$\xymatrixcolsep{2.1pc}
\xy\xymatrix{
0 \ar[r]^-{\gamma_0} &1 \ar[r]^-{\alpha_1} \uloopr{\beta_1} &2
\ar[r]^-{\gamma_1} &3
\ar[r]^-{\alpha_2} \uloopr{\beta_2} &4 \ar[r]^-{\gamma_2} &5
\ar[r]^-{\alpha_3} \uloopr{\beta_3} & \ar@{}[r]|\cdots &
\ar[r]^-{\gamma_{m-1}} &2m-1 \ar[r]^-{\alpha_m} \uloopr{\beta_m} &2m
}\endxy$$
}
\noindent set $p_i= \alpha_i\beta_i$, and $p= p_m\gamma_{m-1}\cdots
p_2\gamma_1p_1\gamma_0$. We define relations $r_1,\dots,r_M\in K\Gamma$ via
$r_i= \sum_{A\in\P} c_i(A)r_i(A)$, where $r_i(A)= s_m\gamma_{m-1}\cdots
s_2\gamma_1s_1\gamma_0$ with
$s_j=\alpha_j$ if
$j\in A$ and $s_j= p_j$ otherwise. The finite dimensional algebra $\la
=K\Gamma/\langle r_1,\dots,r_M,\beta_1^2,\dots,\beta_m^2 \rangle$ clearly
satisfies (N), and that
$V_p\cong V$ is verified as in [6, loc\. cit.]. Since $p$ has no right
subpaths of positive length ending in the vertex 0, 
we obtain $\Phi(k)\ne \Phi(l)$ whenever $k$ and $l$ are
distinct points of
$V_p$ (cf\. [6, Theorem B]). In conclusion, $V_\SS$ is isomorphic to $V$,
where 
$$\SS=
(S_0,S_1,S_1,S_2,S_3,S_3,S_4,\dots,S_{2m-1},S_{2m-1},S_{2m})$$
is the sequence of
the simple modules corresponding to the vertices along $p$, and $\Phi_\SS$ is
a bijection.
\qed\enddemo

\head 4. Algebras for which all uniserial varieties are finite\endhead

Clearly, the condition $|V_p|<\infty$ for all paths $p\in K\Gamma$ forces
$\la$ to have finite uniserial type and thus to satisfy (N).
Actually, the condition `$|V_p|<\infty$ for all $p$' is significantly stronger
than finite uniserial type and is no longer left-right symmetric. Alternately,
algebras with this uniserial behavior can be characterized by the fact that
there is a one-to-one correspondence between the isomorphism types of
uniserial left
modules and their graphs. The description which most readily permits to
recognize such algebras is a tightening of the necessary condition (N) for
finite uniserial type which eliminates one of the two potential
shapes of a mast between vertices spanned by an arrow.

\proclaim{Theorem 4.1} Given any algebra $\la= K\Gamma/I$, the following
conditions are equivalent:

{\rm (1)} For each path $p\in K\Gamma$,  the variety $V_p$ is finite.

{\rm (1')} For each path $p\in K\Gamma$, the variety $V_p$ is either empty or
a singleton.

{\rm (2)} Whenever $\alpha : e\rightarrow e'$ is an arrow in $\Gamma$ and $p
: e\rightarrow e'$ a mast of positive length, the path $p$ is of the form
$p=c'\alpha$, where
$c'$ is an oriented cycle which may be trivial, i.e., $p$ has the form
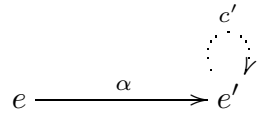
\ignore{
$$\xy\xymatrix{
e \ar[rr]^\alpha &&e' \udotloopr{c'} 
}\endxy$$
}

{\rm (3)} $\la$ has finite uniserial type, and there are no
circular halyards on the mast of any uniserial left $\la$-module.

{\rm (4)} There is a 1--1 correspondence between the isomorphism types and
the graphs of the uniserial left $\la$-modules.

{\rm (5)} The only graphs of uniserial left $\la$-modules are edge
paths.\endproclaim

\remark{Remark} In coordinate-free terms, condition (2) reads as follows: If
there exists a uniserial left $\la$-module $W$ of length 2 with top $S$ and
socle
$S'$, and if $U$ is any uniserial left $\la$-module of length $\ge2$ with top
$S$ and socle
$S'$, then $U/J^2U\cong W$.\endremark

Before we derive the theorem from [6] and the results of the previous
section, we look at examples.
\medskip

\example{Examples 4.2} (a) Whenever $\la= \O/\pi\O$ is the finite dimensional
algebra associated with a tiled classical order $\O$ over a DVR with
uniformizing parameter $\pi$, the equivalent conditions of Theorem 4.1 are
satisfied. Indeed, from Corollary 3.5 we know that $\la$ has finite uniserial
type. In order to guarantee the absence of circular halyards on arbitrary
masts of uniserial left
$\la$-modules, it suffices to recall that no mast runs through an oriented
cycle, since no uniserial module has multiple composition factors. Thus,
condition (3) is satisfied.

(b) Whenever $\la= K\Gamma/I$ for an acyclic quiver $\Gamma$, finite
uniserial type of $\la$ is tantamount with the condition `$|V_p|\le1$ for
all paths $p\in K\Gamma$'. (Compare with Corollary 3.4.)

(c) On the other hand, if $\la= K\Gamma/ \langle \alpha^2\rangle$,
where $\Gamma$ is the quiver
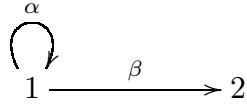
\ignore{
$$\xy\xymatrix{
1 \ar[rr]^\beta \uloopr{\alpha} &&2
}\endxy$$
}
\noindent and $\SS= (S_1,S_1,S_2)$, then $\la$ has finite uniserial type,
whereas
$V_\SS= \AA^1$ is not reduced to a point. Note, however, that the category
mod-$\la$ of finitely generated right $\la$-modules does satisfy the
equivalent conditions of the theorem. In fact, the quickest way to see that
$\lamod$ has finite uniserial type is to note that mod-$\la$ has this
property, and to then apply duality.
\qed\endexample

\demo{Proof of Theorem 4.1} `(1')$\Longrightarrow$(4)' is trivial.

`(4)$\Longrightarrow$(3)'. Finiteness of the uniserial type of $\la$ is
clear. If there were a circular halyard
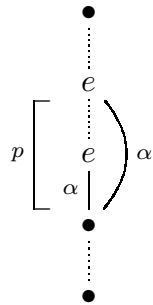
\ignore{
$$\xymatrixrowsep{1pc}
\xy\xymatrix{
\bullet \ar@{{}.{}}[d]\\
e \ar@{{}.{}}[d] \ar@{-}@/^1pc/[dd]^\alpha
\ar@^{|-|}[dd]<-4ex>_p\\
e \ar@{-}[d]_\alpha\\
\bullet \ar@{{}.{}}[d]\\
\bullet}\endxy$$
}
\noindent on some mast -- the subpath $p$ as marked starting in the vertex $e$
say -- then $(\alpha,e)$ would be a slack halyard on the mast $p$ by Lemma
3.1. It is straightforward to check that this would entail the existence of
a uniserial left $\la$-module with mast $p$ having one graph that shows the
halyard
$(\alpha,e)$ and another graph that does not.

The implication `(3)$\Longrightarrow$(5)' follows from Theorem 3.3, and the
implications `(5)$\Longrightarrow$(1')' and `(1')$\Longrightarrow$(1)' are
trivial.

`(2)$\Longrightarrow$(5)'. Let $U$ be a uniserial left $\la$-module with
mast $q$, and let $x\in U$ be any top element. Suppose there is a detour
$(\alpha,u)$ on $q$ such that $\alpha ux\ne 0$, where $\alpha$ is an arrow
$e\rightarrow e'$. Then there is a right subpath $v$ of $q$ longer than $u$
which ends in $e'$, and if we write $v=wu$ for a suitable left subpath $w$
of $v$, then $w : e\rightarrow e'$ is a mast of positive length. That
$(\alpha,u)$ is a detour on $q$ means that $w=w'\beta$, where $\beta$ is an
arrow different from $\alpha$, a situation which is disallowed by condition
(2).

`(3)$\Longrightarrow$(2)'. Assuming (3), we know from Theorem 3.3 that
$\la$ satisfies condition (N). Hence, given any arrow $\alpha : e\rightarrow
e'$,  any mast $p : e\rightarrow e'$ of positive length which fails to
satisfy condition (2) of the theorem has the form
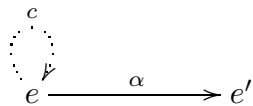
\ignore{
$$\xy\xymatrix{
e \ar[rr]^\alpha \udotloopr{c} &&e'
}\endxy$$
}
\noindent where $c$ is a nontrivial oriented cycle starting in an arrow
$\beta\ne\alpha$. By Lemma 3.1 and condition (N), the detour $(\alpha,e)$ is a
circular halyard on $p$, and condition (3) is violated.

Finally, `(1)$\Longrightarrow$(3)' since given any circular halyard on some
mast, there is a circular halyard of the form
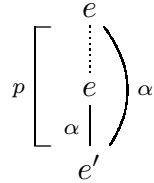
\ignore{
$$\xymatrixrowsep{1.2pc}
\xy\xymatrix{
e \ar@{.}[d] \ar@{-}@/^1pc/[dd]^\alpha
\ar@^{|-|}[dd]<-4ex>_p\\
e \ar@{-}[d]_\alpha\\
e'}\endxy$$
}
\noindent which is necessarily slack by Lemma 3.1. Then Lemma 2.2 tells us
that $V_p$ is infinite.
\qed\enddemo

\head 5. Algebras of finite uniserial type -- geometry
versus representation theory\endhead

Given any path $p\in K\Gamma$, Theorem A of [6] allows us to decide
whether or not $p$ is a mast of a uniserial left $\la$-module. Moreover, in
the positive case, the systems $S_p(X,Y,Z)$ occurring in
Theorem B of [6] provide us with a foolproof computational tool for the
decision whether (up to isomorphism) there are only finitely many uniserial
objects in
$\lamod$ with mast $p$. Indeed, if for $k,l$ we define `$k\sim l$
$\Longleftrightarrow$ the linear system $S_p(k,l,Z)$ is consistent', then the
answer to this last question is `yes' if and only if $V_p$ splits into
finitely many
$\sim$-equivalence classes. So, if we are willing to apply brute computational
force, the sketched results provide a means to settle the question of
whether a given algebra
$\la$ has finite uniserial type, provided that we are in possession of quiver
and relations. However, this approach does not procure a great deal of
qualitative insight into algebras of finite uniserial type, and the
computational ballast involved in the decision is  higher than necessary. 

So, on one hand, we will promote our
theoretical understanding of the uniserial
representation theory at those sequences $\SS$ of simple left $\la$-modules
which do not occur as composition series of infinitely many non-isomorphic
uniserial modules. On the other hand, we will use 
the varieties
$V_\SS$ corresponding to such sequences of simple modules to derive
information about the quiver and relations of $\la$. It appears,
however, that completely pinning down such sequences $\SS$ in
terms of  quiver and relations would not be very meaningful -- neither
from a theoretical nor from a computational point of view -- in that it would
lead into a maze of combinatorial distinctions. Instead, we will approach the
problem from the necessary and sufficient sides, spelling out fairly simple
conditions which `almost' meet. At the end of this section, we present a
much more manageable procedure for recognizing finite uniserial type. In
Sections 6 and 7, we will then proceed to derive convenient equivalent
conditions for finite uniserial type within specific classes of algebras. 

{\it Throughout this section, we assume that
$\la= K\Gamma/I$ satisfies condition (N) of Theorem 3.3, unless we explicitly
forego this condition.}

We first deal with the comparatively easy case where the initial vertex of
the mast considered carries a loop.

\proclaim{Theorem 5.1} Let $\SS= (S(0),\dots,S(l))$ be a sequence of
simple left $\la$-modules with $V_\SS\ne\varnothing$ and
$\Ext(S(0),S(0))\ne 0$. Then the
following conditions are equivalent:

{\rm (1)} There are only finitely many isomorphism types of uniserial left
$\la$-modules with composition series $\SS$.

{\rm (2)} The unique path $p\in K\Gamma$ of length $l$ passing through the
vertices corresponding to the simple modules $S(i)$ in the order prescribed
by $\SS$ is of the form
\ignore{
$$\xy\xymatrix{
e\ar[r]^\gamma \ar@^{|-|}[0,5]<-4ex>_{\gamma^\mu} &e\ar[r]^\gamma &
\ar@{}[r]|\cdots & \ar[r]^\gamma &e\ar[r]^\gamma &e\ar[r]^\alpha
&e_1\ar[r] &e_2\ar[r] & \ar@{}[r]|\cdots & \ar[r] &e_\nu
}\endxy$$
}
\noindent with $\mu,\nu\ge 0$, $\mu+\nu=l$, and $e\notin
\{e_1,\dots, e_\nu\}$; if $\mu\ge 1$ and $\nu\ge 2$, then $e_1\notin
\{e_2,\dots, e_\nu\}$. Moreover, the only slack halyards on $p$ are of the
form
$(\alpha,\gamma^i)$ with $0\le i\le \mu-1$.

If {\rm (2)} is satisfied, and if there are no tight halyards on $p$,
then there exists precisely one uniserial left $\la$-module with
composition series $\SS$.\endproclaim

\remark{Remarks} 1. If the equivalent statements of the theorem hold, the
slack halyards on $p$ are as follows:
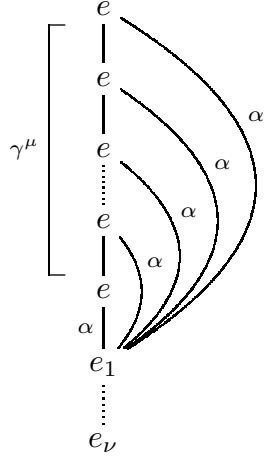
\ignore{
$$\xymatrixrowsep{1pc}
\xy\xymatrix{
e \ar@{-}[d] \ar@^{|-|}[4,0]<-4ex>_{\gamma^\mu}
\ar@{-}@/^4pc/[5,0]^(0.33)\alpha\\
e \ar@{-}[d] \ar@{-}@/^3pc/[4,0]^(0.33)\alpha\\
e \ar@{{}.{}}[d] \ar@{-}@/^2pc/[3,0]^(0.33)\alpha\\
e \ar@{-}[d] \ar@{-}@/^1pc/[2,0]^(0.33)\alpha\\
e \ar@{-}[d]_-\alpha\\
e_1 \ar@{{}.{}}[d]\\
e_\nu
}\endxy$$
}
with every subset of the indicated array of slack halyards
$(\alpha,\gamma^i)$ showing up in one of the graphs of any
uniserial module with composition series $\SS$, for suitable
choices of top elements.

2. A coordinate-free rendering of the second condition of the theorem is as
follows: If
$\mu\ge0$ is maximal with
$S(\mu)\cong S(0)$, then
$S(0)\cong S(1)\cong \cdots \cong S(\mu)$ and, in case $1\le \mu<l$, none of
the simple modules $S(i)$ with $i>\mu+1$ is isomorphic to $S(\mu+1)$;
moreover,
$V_\SS
\cong K^\mu\times V'$ where $V'$ is a finite variety.

Indeed, that this coordinate-free statement implies 5.1(2) in the presence
of condition (N) is clear. For the converse, suppose $V_\SS\subseteq \Bbb A^N$
and, for
$1\le i\le \mu$, let
$X_i$ be the variable corresponding to the slack halyard $(\alpha,
\gamma^{i-1})$. It is straightforward to check the following: Given any
$\mu$-tuple $(l_1,\dots,l_\mu)\in K^\mu$ and any point $k= (k_i)_{i\le N}$
in $V_\SS$, the point $k'= (l_1,\dots, l_\mu, k_{\mu+1},\dots, k_N)$ also
belongs to $V_\SS$. Hence, $V_\SS= K^\mu\times V'$ for a suitable variety
$V'\subseteq \Bbb A^{N-\mu}$ which is finite since, by 5.1(2), the images of
the variables $X_{\mu+1}, \dots, X_N$ in the coordinate ring of $V_\SS$ are
all algebraic over the base field $K$.\endremark

\demo{Proof of Theorem 5.1} `(1)$\Longrightarrow$(2)'. Suppose that (1)
holds. Then the first statement of (2) follows from Theorem 3.3; just
use the fact that each subpath of $p$ is a mast. Next suppose that $\mu\ge1$
and $\nu\ge2$. To see that $e_1\notin \{e_2,\dots,e_\nu\}$, pick $i\in
\{2,\dots,\nu\}$ and consider the right subpath $v_i= \alpha_{i-1}\cdots
\alpha_1 \alpha \gamma^\mu : e\rightarrow e_i$ of $p$ which starts in the
arrow $\gamma \ne \alpha$; here $\alpha_i$ labels the arrow $e_i\rightarrow
e_{i+1}$. Since $e\notin \{e_1,\dots,e_\nu\}$, we have $\alpha_{i-1} \ne
\alpha$, and hence the equality `$e_i=e_1$' is excluded by (N); indeed, (N)
guarantees that each mast $e\rightarrow e_1$ of positive length either
starts in $\alpha$ or ends in $\alpha$.

Now let $(\beta,u)$ be a slack
halyard on $p$; it is necessarily circular by condition (N). So, in case
$u=\gamma^i$ for some
$i$, we must have $\mu\ge1$, $i\le \mu-1$ and $\beta=\alpha$, which puts the
halyard $(\beta,u)$ on the list of admissible slack halyards specified under
(2). Finally, suppose that
$\len(u)>\mu$, i.e., $u=\alpha_m\cdots \alpha_1\alpha \gamma^\mu$ for some
$m\ge0$. Let
$$\beta u \seq X_1v_1 +\dots+ X_sv_s$$
be the substitution equation for $(\beta,u)$, where the $X_i= X_i(\beta,u)$
are independent variables and $v_1,\dots,v_s$ the right subpaths of $p$ of
lengths exceeding that of $u$ and ending in the same vertex $e'\in \{e_1,
\dots, e_\nu\}$ as $\beta$. Since $(\beta,u)$ is slack, one of the above
variables $X_i$ assumes infinitely many $K$-values on $V_p$. By (1), we
infer that there exist distinct points $k,l\in V_p$, with $k_r =k_r(\beta,u)
\ne l_r(\beta,u) =l_r$ for some $r$ say, such that $\Phi_p(k)= \Phi_p(l)$. In
view of [6, Theorem B], we obtain a uniserial left $\la$-module $U$ with
mast
$p$ and top element $y$ such that
$$\sum_{1\le i\le s} (k_i-l_i)v_iy =\sum_{1\le j\le\mu} c_j\beta u\gamma^jy
-\sum_{i,j} k_ic_jv_i
\gamma^jy$$
for certain scalars $c_1,\dots,c_\mu \in K$; indeed,
$\gamma,\dots,\gamma^\mu$ are precisely the right subpaths of $p$ from $e$
to $e$ of positive length. But $u\gamma^jy =\alpha_m\cdots
\alpha_1\alpha \gamma^{\mu+j}y =0$ for $j\ge 1$, and a
fortiori $v_i\gamma^jy=0$, since $u$ is a right subpath of each $v_i$.
Therefore, we conclude that $\sum_i (k_i-l_i)v_iy =0$, which contradicts the
fact that the elements $v_iy$, $1\le i\le s$, are $K$-linearly independent
and thus finishes the proof of `(1)$\Longrightarrow$(2)'.

`(2)$\Longrightarrow$(1)'. It clearly suffices to show that the slack halyards
on $p$ are irrelevant to the isomorphism types of the corresponding
uniserial modules (cf\. Remark 2.3(a)). So let $U$ be a uniserial module
with mast
$p$ and top element $y$. Then $\alpha\gamma^iy =k_i\alpha \gamma^\mu y$ for
$0\le i\le \mu-1$ and suitable scalars $k_i$ because, by (2), each of the
halyards $(\alpha,\gamma^i)$ has a unique variable, namely that
corresponding to the right subpath $v= \alpha\gamma^\mu$ of $p$. We want to
find an alternate top element of the form
$x= y+\sum_{j=1}^\mu c_j\gamma^jy$ of
$U$ such that
$\alpha\gamma^ix =0$ for all $i\le \mu-1$, which amounts to finding scalars
$c_j$ that solve the system
$$k_i\alpha \gamma^\mu y = -\sum_{j\ge 1} c_j\alpha
\gamma^{i+j}y= -c_{\mu-i} \alpha \gamma^\mu y -\sum_{1\le j<
\mu-i} c_jk_{i+j}\alpha \gamma^\mu y$$
$(0\le i\le \mu-1)$ in $U$, since $\gamma^{i+j}y=0$ whenever $i+j>\mu$. The
element $\alpha\gamma^\mu y$ being nonzero, this leads to an equivalent linear
system for the $c_j$ over $K$, the coefficient matrix of which is upper
triangular with entries $-1$ along the main diagonal. The consistency of this
last system completes the proof of the equivalences.

The final statement of the theorem is justified by the proof of the
implication `(2)$\Longrightarrow$(1)'.
\allowlinebreak\qed\enddemo

\example{Example 5.2} Let $\la= K\Gamma/\langle \gamma^4, \alpha_2^3\rangle$
be the monomial relation algebra with quiver
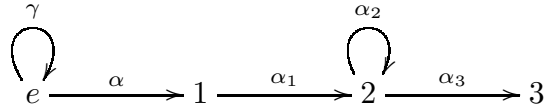
\ignore{
$$\xymatrixcolsep{3.5pc}
\xy\xymatrix{
e \ar[r]^\alpha \uloopr{\gamma} &1 \ar[r]^{\alpha_1} &2
\ar[r]^{\alpha_3} \uloopr{\alpha_2} &3
}\endxy$$
}
\noindent and let $p= \alpha_2\alpha_1\alpha \gamma^3$. It is readily
verified that $p$ is a mast without tight halyards. Theorem 5.1 therefore
guarantees that there is precisely one uniserial left
$\la$-module module
$U$ with mast $p$, up to isomorphism, and that the graphs of $U$ are subgraphs
of
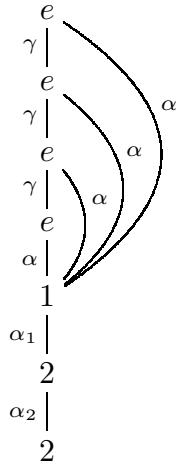
\ignore{
$$\xymatrixrowsep{1pc}
\xy\xymatrix{
e \ar@{-}[d]_\gamma \ar@{-}@/^3pc/[4,0]^(0.35)\alpha\\
e \ar@{-}[d]_\gamma \ar@{-}@/^2pc/[3,0]^(0.35)\alpha\\
e \ar@{-}[d]_\gamma \ar@{-}@/^1pc/[2,0]^(0.35)\alpha\\
e \ar@{-}[d]_\alpha\\
1 \ar@{-}[d]_{\alpha_1}\\
2 \ar@{-}[d]_{\alpha_2}\\
2}\endxy$$
}
\noindent depending on the choice of top element.

On the other hand, there are infinitely many uniserial modules with mast
$q=\alpha_3p$, because $(\alpha_3, \alpha_1\alpha \gamma^3)$ is a slack
halyard on $q$ which is not of the kind specified in part (3) of Theorem
5.1.
\qed\endexample

When we drop the hypothesis that the top simple module have non-trivial
self-ex\-ten\-sions, our problem becomes far more complex. The next two
theorems separately deal with the `necessary' and `sufficient' sides. We
first give a coordinatized formulation which is handy for deciding the
question of whether there are infinitely many uniserial modules with a given
composition series $\SS$, provided that $\la$ is presented by means of quiver
and relations. As usual, given a path
$p$ in
$K\Gamma$ and a detour $(\alpha,u)$ on $p$, the right subpaths of $p$ longer
than $u$ which end in the same vertex as $\alpha$ are denoted by
$v_i(\alpha,u)$, $i\in I(\alpha,u)$; moreover, $X_i(\alpha,u)$ is the
variable going with $v_i(\alpha,u)$.

We will start with a suggestive preview of the following two theorems.

\proclaim{Preview 5.3} Suppose that (N) is satisfied, and let $p\in K\Gamma$
be a mast of positive length, with initial vertex $e$.

If there are only finitely many uniserial modules with mast $p$, then, given
any slack halyard $(\alpha,u)= (\alpha, \beta_\mu\cdots
\beta_1)$ on $p$ with slack variable $X_r(\alpha,u)$, the
corresponding right subpath $v_r(\alpha,u)$ of $p$ looks as follows:
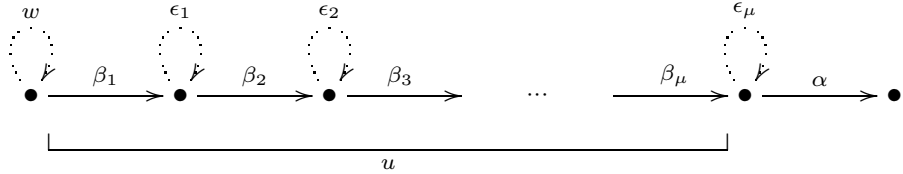
\ignore{
$$\xymatrixcolsep{3pc}
\xy\xymatrix{
\bullet \ar[r]^{\beta_1} \udotloopr{w}
 \ar@^{|-|}[0,5]<-4ex>_u & \bullet
\ar[r]^{\beta_2} \udotloopr{\epsilon_1} & \bullet \ar[r]^{\beta_3}
\udotloopr{\epsilon_2} & \ar@{}[r]|\cdots &
\ar[r]^{\beta_\mu} & \bullet \ar[r]^\alpha
\udotloopr{\epsilon_{\mu}} & \bullet
}\endxy$$
}
\noindent where the oriented cycle $w$ is a right subpath of $p$ of positive
length, whereas the cycles $\epsilon_i$ may be trivial, and $u$
refers to the composition $\beta_\mu \beta_{\mu-1} \cdots \beta_2\beta_1$.
Moreover, if
 $\dim V_p =d$, then the path $p$ passes at least $d+1$ times through its
starting vertex (if we count the initial occurrence of this vertex).

Conversely, suppose that, for each slack halyard $(\alpha,u)$ on $p$ and
each slack variable $X_r(\alpha,u)$, there exists a right subpath $w :
e\rightarrow e$ of $p$ of positive length such that $v_r(\alpha,u)$ has the
form $\alpha uw$, that is,
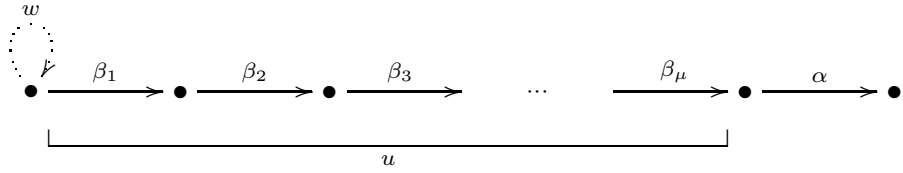
\ignore{
$$\xymatrixcolsep{3pc}
\xy\xymatrix{
\bullet \ar[r]^{\beta_1} \udotloopr{w}
 \ar@^{|-|}[0,5]<-4ex>_u & \bullet
\ar[r]^{\beta_2} & \bullet \ar[r]^{\beta_3}
& \ar@{}[r]|\cdots &
\ar[r]^{\beta_\mu} & \bullet \ar[r]^\alpha
& \bullet
}\endxy$$
}
\noindent Then there are only finitely many uniserial left $\la$-modules
with mast $p$.\endproclaim

Note that the necessary and sufficient conditions above differ only by the
presence or absence of the cycles $\epsilon_i$. In general, these
cycles cannot be deleted from the necessary condition, as Example 5.6(b)
below shows. On the other hand, examples can be constructed to show that our
sufficient conditions cannot be relaxed so as to allow for
cycles $\epsilon_i$, in general.

For a coordinate-free formulation, the following terminology is convenient: a
{\it segment} of a sequence $(a_0,\dots,a_l)$ is a subsequence
$(a_r,a_{r+1},\dots,a_s)$ for indices $r$ and $s$ with $0\le r<s\le l$; we
will speak of a {\it left segment} if $r=0$, of a {\it right segment} if
$s=l$.

As an illustration, we reformulate the second of the above conditions
without relying on coordinates; cf\. end of Section 2. If $\SS
=(S(0),\dots,S(l))$ is a sequence of (isomorphism types of) simple left
$\la$-modules, then there are only finitely many uniserials with composition
series $\SS$, provided that the following condition is satisfied for each
slack halyard
$(S,\rho)$ of
$\SS$ with slack variables
$\sigma_1, \sigma_2,\dots, \sigma_m$. Namely, for each
$i$ between 1 and $m$, the sequence
$(S(0),\dots, S(\rho), S)$ is a
proper {\it right} segment of the sequence $(S(0),\dots, S(\sigma_i))$, up
to coordinatewise isomorphism of the entries.

\proclaim{Theorem 5.4} We continue to assume (N). Let $p= \beta_l
\cdots\beta_1 \in K\Gamma$ be a mast of positive length $l$, consisting of
arrows $\beta_i : e(i-1)\rightarrow e(i)$, such that the set of isomorphism
types of uniserial left $\la$-modules with mast $p$ is finite.
Then: Given any slack halyard $(\alpha,u)= (\alpha, \beta_\mu\cdots
\beta_1)$ on $p$, with a slack variable $X_r(\alpha,u)$ say, the following
is true: 

{\rm({\bf Nec},\,$p$)} There exists a right subpath $w : e(0)\rightarrow e(0)$
of
$p$ of positive length such that
$v_r(\alpha,u)$ is of the form
$$v_r(\alpha,u) =\alpha \epsilon_{\mu}\beta_\mu \epsilon_{\mu-1} \cdots
\epsilon_2\beta_2 \epsilon_1\beta_1w$$
in $K\Gamma$, where each $\epsilon_i$ is an oriented cycle
$e(i)\rightarrow e(i)$ which may be trivial (i.e., $\epsilon_i =e(i)$
is allowed).

Moreover, if $\dim V_p=d$, the path $p$ passes at least $d+1$
times through the vertex $e(0)$, including the start.\endproclaim

\remark{Remark} In
coordinate-free terms, the last part of the theorem says the following:
Suppose that there are only finitely many uniserial modules with composition
series
$\SS= (\la e(0)/Je(0),
\dots,
\la e(l)/Je(l))$. If
$\dim V_\SS =d$, then the simple module $\la e(0)/Je(0)$ occurs with
multiplicity at least $d+1$ in
$\SS$.
\endremark

We prepare for the proof of the theorem with a lemma. Let 
$$D= \{(\alpha,u,i) \mid (\alpha,u)\detour p \text{\ and\ } i\in
I(\alpha,u)\},$$
and recall that $V_p$ is defined as an affine subvariety of $\AA^D$.
Let $w_1,\dots,w_t$ be the distinct right subpaths of $p$ of
positive length of the form
$e(0)\rightarrow e(0)$; then $t$ is clearly the number of recurrences of the
starting vertex $e(0)$ in the sequence of vertices consecutively touched by
the path $p$.

\proclaim {Lemma 5.5} Let $p\in K\Gamma$ be any path (we do not require the
hypotheses of the theorem). Then each fibre of $\Phi_p$ is
contained in a subvariety of $V_p$ of dimension $\le t$.\endproclaim

\demo{Proof} For each $(\alpha,u)\detour p$, we let $I(\alpha,u)$ be a set
of natural numbers, say $I(\alpha,u)= \{1,\dots,\nu(\alpha,u)\}$, such that
$\len v_i(\alpha,u) <\len v_j(\alpha,u)$ whenever $i<j$. Moreover, we equip
the set $D$ of triples with a total order `$\le$' such that $\len
v_i(\alpha,u)< \len v_r(\beta,v)$ implies $(\alpha,u,i) <(\beta,v,r)$.

Let $k^{(0)}$ be a point of $V_p$, and let $\F= \Phi_p^{-1}(\Phi_p(k^{(0)}))$
be the corresponding fibre of
$\Phi_p$. Pick a uniserial module
$U$ in the class 
$\Phi_p(k^{(0)})$, together with a top element $x\in U$ such that $\alpha
ux= \sum k_i^{(0)}(\alpha,u)v_i(\alpha,u)x$ for all $(\alpha,u)\detour p$.
In dependence on $k^{(0)}$, we next define scalars $a_{ij}(\alpha,u)\in K$
via
$$\alpha uw_jx =\sum_{i\in I(\alpha,u)} a_{ij}(\alpha,u)v_i(\alpha,u)x$$
for $(\alpha,u)\detour p$ and $j\le t$. These equations uniquely determine
the $a_{ij}(\alpha,u)$, for the elements $v_i(\alpha,u)x\in U$, $i\in
I(\alpha,u)$, are $K$-linearly independent. Moreover, note that, given $j\le
t$ and $\mu\in I(\alpha,u)$, the path $v_\mu(\alpha,u)w_j$ is strictly
longer than $v_\mu(\alpha,u)$, whence
$$v_\mu(\alpha,u)w_jx =\sum_{i\ge \mu+1} b_{i\mu j}(\alpha,u)v_i(\alpha,u)x$$
for suitable, uniquely determined scalars $b_{i\mu j}(\alpha,u)$; keep in
mind that we have indexed the paths $v_i(\alpha,u)$ so as to reflect their
lengths.

By [6, Section 4], a point $k\in V_p$ belongs to the fibre $\F$ if and
only if there exist scalars $c_1,\dots,c_t$ such that
$$\multline \sum_{i\in I(\alpha,u)} \bigl( k_i(\alpha,u) -k_i^{(0)}(\alpha,u)
\bigr) v_i(\alpha,u)x =\\ \sum_{i\in I(\alpha,u)} \sum_{j\le t}
c_ja_{ij}(\alpha,u) v_i(\alpha,u)x - \sum_{i\in I(\alpha,u)} \sum \Sb \mu\in
I(\alpha,u)\\
\mu\le i-1\endSb \sum_{j\le t} k_\mu(\alpha,u) c_jb_{i\mu j}(\alpha,u)
v_i(\alpha,u)x\endmultline$$
for all $(\alpha,u)\detour p$. This in turn is tantamount with the
solvability of the following $K$-system for $(Z_j)$, $j\le t$: namely,
$$k_i(\alpha,u) -k_i^{(0)}(\alpha,u) =\sum_{j\le t} Z_j \biggl(
a_{ij}(\alpha,u) -\sum \Sb \mu\in I(\alpha,u)\\ \mu\le i-1\endSb
k_\mu(\alpha,u) b_{i\mu j}(\alpha,u) \biggr), \tag{$*_k$}$$
where $(\alpha,u)$ runs through the detours on $p$ and $i$ traces
$I(\alpha,u)$. Given $k\in \F$, we consider, for each triple
$(\alpha,u,i)\in D$, the following portion $M(k,\alpha,u,i)$ of the
coefficient matrix of the system $(*_k)$, namely
$$M(k,\alpha,u,i) =\biggl( a_{rj}(\beta,v) -\sum \Sb \mu\in I(\beta,v)\\
\mu\le r-1\endSb k_\mu(\beta,v) b_{r\mu j}(\beta,v) \biggr) _{(\beta,v,r)\le
(\alpha,u,i);\ j\le t}.$$ 
This matrix has $t$ columns, and rows
indexed by the triples $(\beta,v,r)\le (\alpha,u,i)$. We denote the
rank of the full coefficient matrix of
$(*_k)$ by $m=m(k)$. Clearly, we have $m\le t$. Finally, we define two
sequences
$$\align S(k) &= ((\alpha_1,u_1,i_1),\dots, (\alpha_m,u_m,i_m))\\
T(k) &= (j_1,\dots,j_m), \endalign$$
where the triples $(\alpha_1,u_1,i_1) <(\alpha_2,u_2,i_2)< \cdots
<(\alpha_m,u_m,i_m)$ in $D$ are chosen such that
each $(\alpha_\nu,u_\nu,i_\nu)$ is minimal with respect to the property that
$\rank M(k,\alpha_\nu,u_\nu,i_\nu) =\nu$, and where the $j_\rho$ are
distinct elements in the set $\{1,\dots,t\}$, chosen in such a
fashion that the columns with indices $j_1,\dots,j_\nu$ span the column
space of the matrix $M(k,\alpha_\nu,u_\nu,i_\nu)$ for each $\nu\le m$.

For each pair $(S,T)$ consisting of an ascending sequence $S$ of triples in
$D$ and a sequence $T$ of the same length of elements in $\{1,\dots,t\}$, we
will now construct a closed subvariety
$W(S,T)$ of $V_p$ of dimension $\le t$ which contains all the points $k\in\F$
with $S(k)=S$ and $T(k)=T$. Since $D$ is finite, this will prove the lemma. So
fix a subset $S= \{(\alpha_\nu,u_\nu,i_\nu) \mid 1\le\nu\le m\}$ of $D$ with
$(\alpha_\nu,u_\nu,i_\nu) <(\alpha_{\nu+1},u_{\nu+1},i_{\nu+1})$ for $\nu\le
m-1$, as well as an $m$-tuple $(j_1,\dots,j_m)$ of distinct indices in
$\{1,\dots,t\}$. To find a variety
$W(S,T)$ as claimed, it clearly suffices to show that, for each triple
$(\alpha,u,i)\in D$, there exists a rational function
$f_{(\alpha,u,i)}\in K(\tau_1,\dots,\tau_m)$ which depends only on $k^{(0)}$
and
$S$, $T$, such that $k_i(\alpha,u)= f_{(\alpha,u,i)}(
k_{i_1}(\alpha_1,u_1),\dots, k_{i_m}(\alpha_m,u_m))$ whenever $k\in \F(S,T) :=
\{k'\in\F \mid S(k')=S \text{\ and\ } T(k')=T\}$. We proceed by
induction on the place of $(\alpha,u,i)$ under the ordering on $D$.

First let $(\alpha,u,i)$ be the smallest element of $D$. If $(\alpha,u,i)<
(\alpha_1,u_1,i_1)$, then clearly $k_i(\alpha,u) =k_i^{(0)}(\alpha,u)$ for
all $k\in\F(S)$ by
$(*_k)$ and we let $f_{(\alpha,u,i)}$ be the constant polynomial
$k_i^{(0)}(\alpha,u)$. If $(\alpha,u,i)= (\alpha_1,u_1,i_1)$, set
$f_{(\alpha,u,i)}= \tau_1$. Now suppose that $(\alpha,u,i)$ is not the least
element of $D$. If $(\alpha,u,i)< (\alpha_1,u_1,i_1)$, we again obtain
$k_i(\alpha,u)= k_i^{(0)}(\alpha,u)$ for
all $k\in\F(S)$ and are done. So let $s\ge1$ be the
largest of the numbers $1,\dots,m$ such that $(\alpha,u,i)\ge
(\alpha_s,u_s,i_s)$. Define $f_{(\alpha,u,i)} =\tau_s$ if $(\alpha,u,i)
=(\alpha_s,u_s,i_s)$, and focus on the situation where $(\alpha,u,i)
>(\alpha_s,u_s,i_s)$. Then
$\rank M(k,\alpha,u,i)=s$ for each $k\in\F(S)$, by construction, and the row
$R(k,\alpha,u,i)$ with label
$(\alpha,u,i)$ of the matrix
$M(k,\alpha,u,i)$ is a $K$-linear combination of the rows
$R(k,\alpha_1,u_1,i_1),\dots, R(k,\alpha_s,u_s,i_s)$ of that matrix, say
$$R(k,\alpha,u,i)= \sum_{\nu\le s} \rho_\nu(k) R(k,\alpha_\nu,u_\nu,i_\nu)$$
with $\rho_1(k),\dots,\rho_s(k)\in K$.

Let us scrutinize the linear system
$$R(k,\alpha,u,i) =\sum_{\nu\le s} Y_\nu R(k,\alpha_\nu,u_\nu,i_\nu)$$
for the $\rho_\nu(k)$. Because the $s\times t$ matrix with rows
$R(k,\alpha_1,u_1,i_1),\dots, R(k,\alpha_s,u_s,i_s)$ has rank $s$, it has a
unique solution $(y_1,\dots,y_s)= (\rho_1(k),\dots,\rho_s(k))$. Since,
moreover, the column space of the coefficient matrix is spanned by the
columns labelled $j_1,\dots,j_s$, Cramer's Rule allows us to express the
$y_\bullet$'s as rational functions in the $j_1,\dots,j_s$-entries of the
rows
$R(k,\alpha,u,i)$ and
$R(k,\alpha_\nu,u_\nu,i_\nu)$ with $\nu\le s$. A typical coefficient of the
former row is of the form
$$a_{ij}(\alpha,u) -\sum \Sb \mu\in I(\alpha,u)\\ \mu\le i-1\endSb
f_{(\alpha,u,\mu)}( k_{i_1}(\alpha_1,u_1),\dots, k_{i_m}(\alpha_m,u_m))
b_{i\mu j}(\alpha,u)$$
by the induction hypothesis, where the
$f_{(\alpha,u,\mu)}(\tau_1,\dots,\tau_m)$ are rational functions only
depending on $k^{(0)}$ and $S$, $T$. Indeed, any triple
$(\alpha,u,\mu)$ with
$\mu< i$ is smaller than $(\alpha,u,i)$ in our ordering. In other words, this
coefficient is of a form 
$$g(k_{i_1}(\alpha_1,u_1),\dots,
k_{i_m}(\alpha_m,u_m))$$
for a rational function $g$ which depends only on
$k^{(0)}$ and $S$, $T$. Similarly, the induction hypothesis permits us to
write all coefficients of the rows $R(k,\alpha_\nu,u_\nu,i_\nu)$ with $\nu\le
s$ as rational functions of the required type, and hence each $\rho_\nu(k)$ is
of the form 
$$h(k_{i_1}(\alpha_1,u_1),\dots, k_{i_m}(\alpha_m,u_m))$$
for some
rational function $h\in K(\tau_1,\dots,\tau_m)$ which does not depend on $k$.
But in light of the fact that each $k\in \F(S,T)$ satisfies $(*_k)$, we see
that
$$k_i(\alpha,u) =k_i^{(0)}(\alpha,u) +\sum_{\nu\le s} \rho_\nu(k) \bigl(
k_{i_\nu}(\alpha_\nu,u_\nu) -k^{(0)}_{i_\nu}(\alpha_\nu,u_\nu) \bigr),$$
which completes the induction. \qed\enddemo

\demo{Proof of Theorem 5.4} Let $v_i =v_i(\alpha,u)$, $i\in I(\alpha,u)$, be
the right subpaths of $p$ which are longer than $u$ and end in the same vertex
$e'$ as
$\alpha$. Suppose that the variable $X_r= X_r(\alpha,u)$ in the corresponding
substitution equation
$$\alpha u\seq \sum_{i\in I(\alpha,u)} X_iv_i$$
is slack.  By hypothesis, there exist points $k,k'\in V_p$ with $k_r=
k_r(\alpha,u) \ne k'_r(\alpha,u) =k'_r$ such that $\Phi_p(k)= \Phi_p(k')$. We
apply the method of [6, Theorem B] to obtain a uniserial module $U$ in
this isomorphism class, with top element $y$ say, such that $\alpha uy
=\sum_{i\in I(\alpha,u)} k_iv_iy$, and
$$\sum_{i\in I(\alpha,u)} (k'_i-k_i)v_iy =\sum_{j\le t} c_j\alpha uw_jy
-\sum_{i\in I(\alpha,u)}
\sum_{j\le t} k'_ic_jv_iw_jy,$$
where $c_j\in K$ and the $w_j$ are right subpaths $e(0)\rightarrow e(0)$ of
$p$ of positive length. Each of the summands involved belongs to the
$K$-space $e'J^{\len(u)+1}U$ which has basis $v_iy$, $i\in I(\alpha,u)$. In
view of `$k'_r-k_r\ne 0$', we infer that a nonzero term $kv_ry$ arises on the
right, either in the expansion of some term $\alpha uw_jy$ or that of
some term
$v_iw_jy$. By [6, Theorem A], this implies that either $\alpha uw_j$ or
$v_iw_j$ is a route on $v_r$. This, in turn, forces $v_r$ to be strictly
longer than $w_j$ in both cases, say $v_r= vw_j$ in $K\Gamma$, where $v :
e(0)\rightarrow e'$ is a non-trivial left subpath of $v_r$. In the first
case,
$\alpha u$ is a route on $v$, in the second $v_i$ is a route on $v$. But
$\alpha u$ being a route on each $v_m$, this shows that $\alpha u=
\alpha\beta_\mu
\cdots \beta_1$ is a route on $v$ in both situations.

In particular, $\beta_1$ is a route on $v$, and if we write
$v=v'q_1$ in $K\Gamma$, where $q_1$ is the shortest right subpath of $v$ on
which $\beta_1$ is a route, then $q_1$ ends in $e(1)$ and $\alpha
\beta_\mu \cdots\beta_2$ is a route on $v'$. Observe that $q_1 :
e(0)\rightarrow e(1)$, being a subpath of $p$, is a mast of positive length,
and therefore Condition (N) forces $q_1$ to be of one of the forms
$\beta_1\epsilon_0$, where $\epsilon_0$ is a (possibly trivial) cycle
$e(0)\rightarrow e(0)$, or $\epsilon_1\beta_1$, where $\epsilon_1$ is a
nontrivial cycle $e(1)\rightarrow e(1)$. But our minimal choice of $q_1$
rules out the second possibility, and so $q_1= \beta_1\epsilon_0$. Continuing
inductively, we obtain a factorization $v= v''\beta_\mu\epsilon_{\mu-1}
\cdots\beta_1\epsilon_0$, where each $\epsilon_i$ is a cycle
$e(i)\rightarrow e(i)$, and $(\alpha,e(\mu))$ is a route on the path $v''
: e(\mu) \rightarrow e'$. In particular, $v''$ has positive length, and it
suffices to show that $v$ ends in the arrow $\alpha : e(\mu) \rightarrow
e'$ to establish our claim with $w= \epsilon_0w_j$. To see this, note that
$v_r=vw_j$ contains $u$ as a proper right subpath, say $v_r=u'u$ in
$K\Gamma$; here $u' : e(\mu)\rightarrow e'$ is a mast of positive length
which does not start in $\alpha$ since $(\alpha,u)$ is a detour on $p$ and
hence on $v_r$. Invoking Condition (N) once more, we conclude that $u'$ ends
in $\alpha$. So do $v_r$ and $v$, as a consequence. This proves the first
part of the theorem.

For the last part, let $t$ again be the number of right subpaths
$e(0)\rightarrow e(0)$ of
$p$ of positive length. Since, due to our hypothesis on $p$, the variety
$V_p$ is covered by finitely many fibres of $\Phi_p$,  Lemma 5.5
guarantees that $V_p$ is covered by a finite number of closed subvarieties
of dimension $\le t$. This proves $\dim V_p\le t$ as desired.
\qed\enddemo

\example{Remark and Example 5.6} Let us return to the graphical rendering of
the first part of this theorem given in the Preview. If $\la$ is of finite
uniserial type and thus, in particular, satisfies (N), then any mast $p$ with
a slack halyard of the form
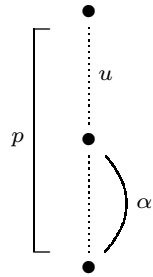
\ignore{
$$\xymatrixrowsep{1pc}
\xy\xymatrix{
\bullet \ar@{{}.{}}[dd]^u \ar@^{|-|}[4,0]<-4ex>_p\\
\\
\bullet \ar@{{}.{}}[dd] \ar@{-}@/^1pc/[2,0]^\alpha\\
\\
\bullet}\endxy$$
}
\noindent where $u= \beta_\mu \cdots\beta_1$, has a right subpath which looks
as follows:
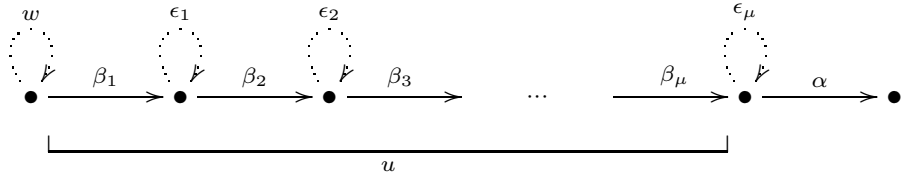
\ignore{
$$\xymatrixcolsep{3pc}
\xy\xymatrix{
\bullet \ar[r]^{\beta_1} \udotloopr{w}
 \ar@^{|-|}[0,5]<-4ex>_u & \bullet
\ar[r]^{\beta_2} \udotloopr{\epsilon_1} & \bullet \ar[r]^{\beta_3}
\udotloopr{\epsilon_2} & \ar@{}[r]|\cdots &
\ar[r]^{\beta_\mu} & \bullet \ar[r]^\alpha
\udotloopr{\epsilon_{\mu}} & \bullet
}\endxy$$
}
\noindent Here the cycles $\epsilon_i$ may be trivial, whereas $w$ is
not. Note that either $u$ is a right subpath of $w$ or vice versa, since both
$u$ and
$w$ are right subpaths of $p$, which shows that $p$ includes at least two
distinct nontrivial oriented cycles, possibly overlapping. This imposes
a high degree of repetitiveness on masts with slack halyards in the case of
finite uniserial type.

As for the converse: Roughly speaking, the
next theorem will guarantee that $\la$ has finite uniserial type in case the
$\epsilon_i$ are trivial for all slack halyards $(\alpha,u)$ and masts $p$ as
above. In case no trivial choice of certain
$\epsilon_i$ is possible, these cycles provide further circular halyards --
usually slack on suitable subpaths $q$ of $p$; whether $q$ supports
infinitely many uniserial modules, can then in turn be tested in the light of
Theorem 5.4. This makes it far less involved, computationally, to decide
the question of finite uniserial type than that of whether a specified mast
$p$ supports infinitely many uniserial modules. Another fact which often
simplifies the decision process as to whether $\la$ has finite uniserial
type is the obvious left-right symmetry of this condition.

(a) For instance, let $\la= K\Gamma/I$, where $\Gamma$ is the quiver
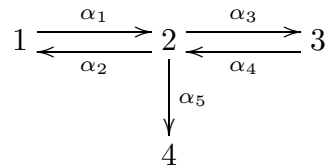
\ignore{
$$\xymatrixcolsep{3pc}
\xy\xymatrix{
1 \ar[r]<0.75ex>^{\alpha_1} &2 \ar[l]<0.75ex>^{\alpha_2}
\ar[r]<0.75ex>^{\alpha_3} \ar[d]^{\alpha_5} &3
\ar[l]<0.75ex>^{\alpha_4}\\
 &4
 }\endxy$$
}
\noindent and $I$ the ideal generated by $\alpha_2\alpha_1\alpha_2$,
$\alpha_3\alpha_4$, $\alpha_2\alpha_4$, and $\alpha_5\alpha_1\alpha_2
-\alpha_5\alpha_4\alpha_3\alpha_1\alpha_2$. Then $p=
\alpha_5\alpha_4\alpha_3\alpha_1\alpha_2\alpha_1$ has two slack halyards
$(\alpha_5,\alpha_1)$ and $(\alpha_3,\alpha_1)$
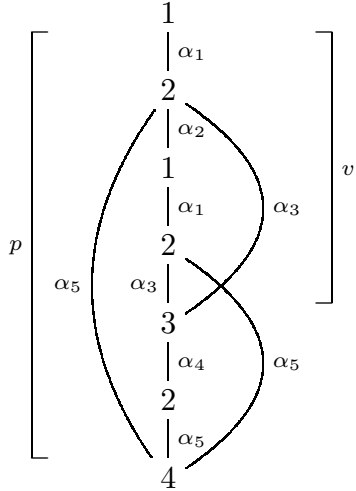
\ignore{
$$\xymatrixrowsep{1pc}
\xy\xymatrix{
1 \ar@{-}[d]^{\alpha_1} \ar@^{|-|}[6,0]<-10ex>_p
\ar@_{|-|}[4,0]<12ex>^v\\
2 \ar@{-}[d]^{\alpha_2} \ar@{-}@/_2pc/[5,0]_{\alpha_5}
\ar@{-}@/^2.5pc/[3,0]^{\alpha_3}\\
1 \ar@{-}[d]^{\alpha_1}\\
2 \ar@{-}[d]_{\alpha_3} \ar@{-}@/^2.5pc/[3,0]^{\alpha_5}\\
3 \ar@{-}[d]^{\alpha_4}\\
2 \ar@{-}[d]^{\alpha_5}\\
4}\endxy$$
}
\noindent and one tight halyard, namely $(\alpha_5,
\alpha_1\alpha_2\alpha_1)$, tightened by the last of the relations. The mast
$p$ satisfies the conditions of Theorem 5.4, for $p$ and $v$ have the
forms
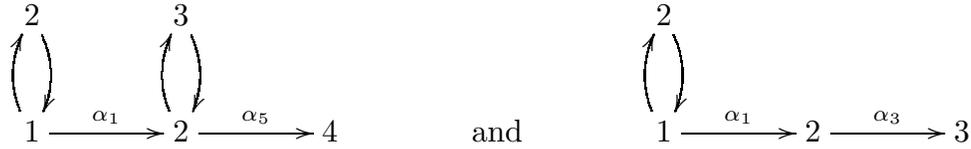
\ignore{
$$\xymatrixcolsep{3pc}
\xy\xymatrix{
2 \ar@/^0.5pc/[d] & 3 \ar@/^0.5pc/[d] &&& 2 \ar@/^0.5pc/[d] \\
1 \ar@/^0.5pc/[u] \ar[r]^{\alpha_1} & 2 \ar@/^0.5pc/[u]
\ar[r]^{\alpha_5} & 4 & \text{and} & 1 \ar@/^0.5pc/[u]
\ar[r]^{\alpha_1} & 2 \ar[r]^{\alpha_3} & 3
}\endxy$$
}
\noindent respectively, and so no immediate information about the uniserial
type of $\la$ is obtained. However, Theorem 5.4 applies to the `dual' mast
$p'$ of the corresponding uniserial {\it right} $\la$-modules $D(U)$, where
$U$ runs through the uniserial left $\la$-modules with mast $p$ and $D :
\lamod\rightarrow \text{mod-}\la$ denotes the standard duality. Observe that
the halyard $(\alpha_1,\alpha_5)$ on $p'$ is slack, the variable $X_2$ in the
substitution equation $\alpha_1\alpha_5 \seq
X_1\alpha_1\alpha_3\alpha_4\alpha_5 +X_2p'$ being slack since $V_{p'}=
V(X_1-1)\cong \AA^1$. Theorem 5.4 shows that there are infinitely many
non-isomorphic uniserial right $\la$-modules with mast $p'$, because
$p'$ has no nontrivial subpaths $e_4\rightarrow e_4$. Therefore, there are
of course infinitely many non-isomorphic uniserial left $\la$-modules with
mast $p$. A fortiori, $\la$ has infinite uniserial type. This last fact can
be recognized without any computational effort: The halyard
\ignore{
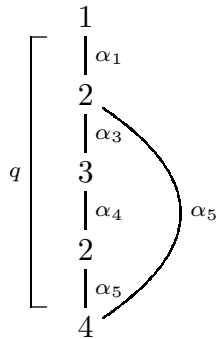
$$\xymatrixrowsep{1pc}
\xy\xymatrix{
1 \ar@{-}[d]^{\alpha_1} \ar@^{|-|}[4,0]<-4ex>_q\\
2 \ar@{-}[d]^{\alpha_3} \ar@{-}@/^2.5pc/[3,0]^{\alpha_5}\\
3 \ar@{-}[d]^{\alpha_4}\\
2 \ar@{-}[d]^{\alpha_5}\\
4}\endxy$$
}
\noindent on the mast $q$ is clearly slack and thus, by Theorem 5.4,
gives rise to infinitely many isomorphism classes of uniserial modules,
since $q$ has no right subpath of positive length which ends in 1.

It is worthwhile to observe that there is only one uniserial module with mast
$\alpha_5\alpha_4\alpha_3\alpha_1\alpha_2$ by Theorem 3.3(IId), since the
halyard $(\alpha_5, \alpha_1\alpha_2)$ is still tight, and the two slack
halyards $(\alpha_5,e_2)$ and $(\alpha_3, e_2)$ emerge from the top vertex.
This illustrates how one may alternately gain and lose the property that
there are only finitely many uniserial modules with a given mast as one lets
the pertinent mast grow (or shrink).

(b) If we strengthen the last of the above relations slightly to tighten the
halyard on the mast $q$, namely, if we let $\Gamma$ be as above and $\la'=
K\Gamma/I'$, where $I'$ is generated by $\alpha_2\alpha_1\alpha_2$,
$\alpha_3\alpha_4$, $\alpha_2\alpha_4$, and $\alpha_5\alpha_1
-\alpha_5\alpha_4\alpha_3\alpha_1$, then $\la'$ has finite uniserial type.
The only `critical' mast being $p$, we use the system $S_p$ of [6,
Theorem B] to ascertain that, this time, there is a single uniserial object
in $\la'$-mod with mast $p$, up to isomorphism. This example also shows that
no trivial choice of the oriented cycles $\epsilon_i$ occurring in the
conclusion of Theorem 5.4 needs to be possible. Indeed, for the slack
halyard $(\alpha,u)= (\alpha_5,\alpha_1)$ on $p$, the only choices of $w$
and $\epsilon_1$ as prescribed by the theorem are $w=\alpha_2\alpha_1$ and
$\epsilon_1= \alpha_4\alpha_3$.
\qed\endexample 

For the following theorem we need not postulate condition (N) in full
strength as an a priori hypothesis. It actually suffices to require that
$\Gamma$ be without double arrows. Observe that, for any
detour $(\alpha,u)$ on $p$, this entails $\len(\alpha u)<
\len(v_i(\alpha,u))$.

\proclaim{Theorem 5.7} Let
$p\in K\Gamma$ be a mast of positive length starting in the vertex $e$ and
satisfying the following condition:

{\rm({\bf Suf},\,$p$)} Whenever $(\alpha,u)$ is a slack halyard on $p$ and
$X_r(\alpha,u)$ a slack variable, there exists a right subpath $w_r(\alpha,u)
: e\rightarrow e$ of
$p$ such that $v_r(\alpha,u) =\alpha uw_r(\alpha,u)$ in
$K\Gamma$.

Then there are only finitely many uniserial left $\la$-modules with mast
$p$, up to isomorphism. \endproclaim

\remark{Remark} For a coordinate-free phrasing of this condition, see the
comments following Preview 5.3.\endremark

\demo{Proof} Set 
$$\align \S &=\{(\alpha,u,i) \mid (\alpha,u) \text{\ is a slack
halyard on\ } p \text{\ and\ } X_i(\alpha,u) \text{\ is a slack variable}\}\\
\T &=\{(\alpha,u,i) \mid (\alpha,u) \text{\ is a slack
halyard on\ } p \text{\ and\ } X_i(\alpha,u) \text{\ is a tight
variable}\},\endalign$$
 and consider the map
$$\align \Psi\ :\ \S\ &\rightarrow \{ \text{right subpaths\ }
e\rightarrow e \text{\ of\ } p \text{\ of positive length} \}\\
(\alpha,u,i) &\mapsto w_i(\alpha,u)\endalign$$
with $v_i(\alpha,u) =\alpha uw_i(\alpha,u)$ as in ({\bf Suf},\,$p$). This
function is well-defined: for each triple $(\alpha,u,i)\in \S$, the path
$w_i(\alpha,u)$ has positive length, because $\alpha u$ is not a right
subpath of $p$, and clearly
$w_i(\alpha,u)$ is uniquely determined by the triple
$(\alpha,u,i)$. This implies in particular that $\len(v_i(\alpha,u))\ge
\len(u)+2$ whenever the triple $(\alpha,u,i)\in\S$.

In a first step, we will prove that
$\Psi$ is an injection. More precisely, we will show by induction on $\len(q)$
that, for each non-trivial right subpath
$q$ of $p$, the restriction of $\Psi$ to the set
$$\S_q =\{ (\alpha,u,i)\in \S \mid \len(v_i(\alpha,u)) \le \len(q)\}$$
is injective.

The case `$\len(q)=1$' is trivial, for the paths $v_i(\alpha,u)$ all have
length at least 2. So let us assume that
$\len(q)\ge 2$. In case $q= v_i(\alpha,u)$ with
$(\alpha,u,i)\in\S$, it will be convenient to assume that $i=q$, that is, 
$q= v_q(\alpha,u)$; moreover, we will assume that $i\ne q$ whenever $q\ne
v_i(\alpha,u)$. 
Using this convention, we obtain: If $q$ ends in the arrow
$\gamma$ and
$(\alpha,u)$ is a detour on $p$, then, clearly, $(\alpha,u,q)\in\S$ only if
$\alpha =\gamma$ and $\len(u)<
\len(q)$ by hypothesis and the above remark. 

Distinct elements of the form $(\gamma,u_1,q),\ (\gamma,u_2,q)\in \S_q$
obviously give rise to distinct paths $w_q(\gamma,u_1)$ and $w_q(\gamma,u_2)$.
To establish injectivity of $\Psi$ on $\S_q$, it thus suffices -- in view of
the induction hypothesis -- to show that the intersection
$$\{w_q(\gamma,u) \mid (\gamma,u,q)\in \S_q\} \cap \{w_i(\beta,r) \mid
(\beta,r,i)\in \S_q,\ i\ne q\}$$
is empty. Suppose, to the contrary, that $w_q(\gamma,u) =w_i(\beta,r) =:w$
for some $(\gamma,u,q)$ and $(\beta,r,i)$ in $\S_q$ with $i\ne q$. From the
equality $\gamma uw=q$ and the fact that $\beta rw$ is a proper right
subpath of $q$, we glean that $\beta r$ is a right subpath of $u$. Write
$u=u'\beta r$ in $K\Gamma$, where $u'$ is a path
which may be trivial. Then $p=p'\beta r$ for a suitable left subpath $p'$
of $p$, which is incompatible with the fact that $(\beta,r)$ is a detour
on
$p$. Thus the first claim is established.

To prove the finiteness of the number of isomorphism classes of uniserial
left $\la$-modules with mast $p$, it clearly suffices to procure, for each
such uniserial module $U$, a top element $x\in U$ with the property that
$\alpha ux\in \sum_{(\alpha,u,i)\in\T} Kv_i(\alpha,u)x$ for every slack
halyard
$(\alpha,u)$ on
$p$ (see Remark 2.3(a)). So let
$U$ be uniserial with mast
$p$.

Due to the injectivity of $\Psi$, we can introduce a total order on the set
$\S$ by defining
$$(\alpha,u,i) <(\beta,v,j)\quad \Longleftrightarrow\quad
\len(w_i(\alpha,u)) <\len(w_j(\beta,v)).$$
In particular, this definition entails that, given distinct indices $i,j\in
I(\alpha,u)$ with $(\alpha,u,i)$ and $(\alpha,u,j)$ in $\S$, we have
$(\alpha,u,i) <(\alpha,u,j)$ if and only if
$v_i(\alpha,u)$ is shorter than $v_j(\alpha,u)$.

We will show that, for any initial segment $\A\subseteq \S$, there exists a
top element $x_\A$ of $U$ such that
$$\alpha ux_\A \in
\sum_{(\alpha,u,i)\in (\S\cup\T)\setminus\A} Kv_i(\alpha,u)x_\A$$
for all slack halyards $(\alpha,u)$ on $p$. Once this is established, the
case $\A=\S$ will finish the proof, since the top element $x=x_\S$ of $U$
will then have the property required above. We will proceed by induction on
$|\A|$.

If $|\A|=0$, i.e., if $\A=\varnothing$, any choice of a top element $x_\A$
of $U$ is as desired.

 So suppose that $|\A|\ge 1$, and let
$(\epsilon,s,\nu)$ be the largest element of $\A$. Applying the induction
hypothesis to the initial segment $\B= \A\setminus \{(\epsilon,s,\nu)\}$
yields a top element $x_\B$ of $U$ such that,
for any slack halyard $(\alpha,u)$ on $p$, we have
$$\align \alpha ux_\B &=\sum_{(\alpha,u,i)\in (\S\cup\T)\setminus \B}
l_i(\alpha,u)v_i(\alpha,u)x_\B\\
 &=\sum \Sb (\alpha,u,i)\in\S\\ (\alpha,u,i)\ge (\epsilon,s,\nu)\endSb
l_i(\alpha,u) v_i(\alpha,u)x_\B + \sum_{(\alpha,u,i)\in\T} l_i(\alpha,u)
v_i(\alpha,u)x_\B\endalign$$   
for
suitable scalars
$l_i(\alpha,u)$. We now define
$x_\A= x_\B -l_\nu(\epsilon,s)w_\nu(\epsilon,s)x_\B$ and check that our
requirements are met. First, we note that
$$\align \epsilon sx_\A &= -l_\nu(\epsilon,s) \epsilon sw_\nu(\epsilon,s)x_\B
+ \sum_{(\epsilon,s,i)\in (\S\cup\T)\setminus \B} l_i(\epsilon,s)
v_i(\epsilon,s)x_\B\\
 &= \sum_{(\epsilon,s,i)\in (\S\cup\T)\setminus\A}
l_i(\epsilon,s) v_i(\epsilon,s)x_\B\\
 &= \sum_{(\epsilon,s,i)\in (\S\cup\T)\setminus\A} l_i(\epsilon,s)
v_i(\epsilon,s)(x_\A +l_\nu(\epsilon,s)w_\nu(\epsilon,s)x_\B).\endalign$$
For each triple $(\epsilon,s,i)\in (\S\cup\T)\setminus\A$, the path
$v_i(\epsilon,s)w_\nu(\epsilon,s)$ is strictly longer than $\epsilon
sw_\nu(\epsilon,s) =v_\nu(\epsilon,s)$ and ends in the same vertex as
$\epsilon$. Therefore  each of the corresponding elements
$v_i(\epsilon,s)w_\nu(\epsilon,s)x_\B$ is a $K$-linear
combination of elements $v_h(\epsilon,s)x_\A$ where $v_h(\epsilon,s)$ is
longer than $v_\nu(\epsilon,s)$, which entails that the triple
$(\epsilon,s,h)$ either belongs to $\T$ or else is strictly larger than
$(\epsilon,s,\nu)$. This shows that
$v_i(\epsilon,s)w_\nu(\epsilon,s)x_\B$ belongs to the
space $\sum_{(\epsilon,s,i)\in (\S\cup\T)\setminus\A} Kv_i(\epsilon,s)x_\A$,
and hence $\epsilon sx_\A \in
\sum_{(\epsilon,s,i)\in (\S\cup\T)\setminus\A} Kv_i(\epsilon,s)x_\A$.

 Finally, we focus on a slack
halyard
$(\alpha,u)$ different from
$(\epsilon,s)$. Note that, for any triple $(\alpha,u,j)\in\S$ with
$(\alpha,u,j)\ge (\epsilon,s,\nu)$, we then have $(\alpha,u,j)>
(\epsilon,s,\nu)$. For any triple $(\alpha,u,j)\in\S$ with $(\alpha,u,j)<
(\epsilon,s,\nu)$, on the other hand, we obtain
$$\len(\alpha uw_\nu(\epsilon,s)) >\len(\alpha uw_j(\alpha,u))
=\len(v_j(\alpha,u))$$
by the definition of our order. Since obviously the path $\alpha
uw_\nu(\epsilon,s)$ ends in the same vertex as $\alpha$, we infer that
$\alpha uw_\nu(\epsilon,s)x_\A$ is a $K$-linear combination of elements
$v_i(\alpha,u)x_\A$ with either $(\alpha,u,i)> (\epsilon,s,\nu)$ or
$(\alpha,u,i)\in\T$. But this means $(\alpha,u,i)\notin \A$.
Consequently, the element
$$\multline \alpha u x_\A =\\
 -l_\nu(\epsilon,s)\alpha uw_\nu(\epsilon,s)x_\A + 
\sum\Sb (\alpha,u,i)\in\S\\ (\alpha,u,i)\ge (\epsilon,s,\nu)\endSb
l_i(\alpha,u)v_i(\alpha,u)x_\A +\sum_{(\alpha,u,i)\in\T}
l_i(\alpha,u)v_i(\alpha,u)x_\A\endmultline$$
belongs to $\sum_{(\alpha,u,i)\in (\S\cup\T)\setminus\A}
Kv_i(\alpha,u)x_\A$. This completes the induction and the proof of the
theorem.
\qed\enddemo

\proclaim{Corollary 5.8 to the proof of Theorem 5.7} Let $p=pe\in K\Gamma$
be a mast of positive length, and suppose that, for \underbar{each} halyard
$(\alpha,u)$ on $p$ and \underbar{each} $i\in I(\alpha,u)$, there is a cycle
$w_i(\alpha,u) : e\rightarrow e$ with $v_i(\alpha,u) =\alpha uw_i(\alpha,u)$
in $K\Gamma$. Then there exists precisely one uniserial left $\la$-module
with mast $p$. Moreover, $V_p\cong \AA^m$, where $m$ is the sum of the
cardinalities $|I(\alpha,u)|$ as $(\alpha,u)$ runs through the halyards on
$p$. \qed\endproclaim

In the next section, we will see that a monomial relation algebra $\la=
K\Gamma/I$ has finite uniserial type if and only if all masts $p$ of
uniserial $\la$-modules (i.e., all paths in $K\Gamma
\setminus I$ in that case) satisfy
$(\bold{Suf},\,p)$. The reader who fills in the details of the next example
will realize how easy it is to apply Theorem 5.7 in the monomial situation.

\example{Example 5.9} Let $\Gamma$ be the quiver
\ignore{
$$\xymatrixcolsep{3pc}
\xy\xymatrix{
2 \ar[r]<0.5ex>^{\alpha_1} &1 \ar[l]<0.5ex>^{\alpha_2}
\ar[r]<0.5ex>^{\alpha_3} &3 \ar[l]<0.5ex>^{\alpha_4}
\ar[r]^{\alpha_5} &4
}\endxy$$
}
\noindent Then the algebra $\la= K\Gamma/ \langle \alpha_2\alpha_1,\
\alpha_3\alpha_4,\ \text{all paths of length\ } \ge9\rangle$ has finite
uniserial type by Theorem 5.7; indeed, it can easily be verified that
each mast $p\in K\Gamma$ of positive length satisfies $(\bold{Suf},\,p)$. In
fact, the stronger hypothesis of Corollary 5.8 is satisfied, whence each
mast $p$ supports exactly one uniserial left $\la$-module, up to isomorphism.
\qed\endexample

\noindent{\bf Algorithm 5.10} for deciding whether a given algebra $\la=
K\Gamma/I$ has finite uniserial type.

By the preceding theory, $\la$ has finite uniserial type if and only if each
of the steps marked with a bullet yields a positive finding.

$\bullet$ Check whether the quiver $\Gamma$ is without double
arrows.

For each path $p$ of length $\ge 2$ in $K\Gamma$, find a set of
polynomials defining the variety $V_p$. Next list those $p$'s for which
$V_p\ne \varnothing$; these are precisely the masts of the uniserial left
$\la$-modules of composition length at least 3.

$\bullet$ Check whether condition (N) of Theorem 3.3 is satisfied.
(Once we have listed the masts, this can be done by mere inspection of
$\Gamma$.)

For the masts $p$ of length $\ge2$, identify the transcendental variables in
the coordinate ring $A_p$.

$\bullet$ For each such mast $p$ and each transcendental variable
of $A_p$, check whether condition $(\bold{Nec},\ p)$ of Theorem 5.4 is
satisfied. (In case the starting vertex of $p$ carries a loop in $\Gamma$,
Theorem 5.1 provides a shortcut.)

$\bullet$ For those masts $p$ which fail to satisfy condition
$(\bold{Suf},\ p)$ of Theorem 5.7, set up the system of equations
$S_p(X,Y,Z)$ of [6, Section 4], and check whether $V_p$ splits into finitely
many equivalence classes under the relation
$$k\sim l \qquad \Longleftrightarrow \qquad \text{the linear system\ }
S_p(k,l,Z) \text{\ is consistent.}$$
\medskip

Finally, we formulate a conjecture which, if confirmed, would complete our
understanding of the interplay between the geometry and representation
theory of uniserial modules. In case the first part of the conjecture can be
proved, the second will follow.

\example{Conjecture 5.11} If $\la$ has finite uniserial type, then 

$\bullet$ All of the nonempty varieties $V_\SS$, where $\SS$ runs through the
sequences of simple left $\la$-modules, are linear, i.e., $V_\SS\cong \AA^m$
for some integer $m$ depending on $\SS$.

$\bullet$ For each sequence $\SS$ of simple left $\la$-modules, there exists
at most one uniserial object in $\lamod$ with composition series
$\SS$.\endexample

\head 6. Monomial relation algebras of finite uniserial type\endhead

In this section, $\la=K\Gamma/I$ will be a monomial relation algebra over an
infinite field $K$; we will assume $I$ to be generated by paths. Moreover,
given a sequence
$(e(0),e(1),\dots,e(l))$ of vertices in $\Gamma$, a {\it path through this
sequence} will be a path $p=
\alpha_l\cdots\alpha_1$ of length $l$, where each $\alpha_i$ is an arrow
$e(i-1)\rightarrow e(i)$.

\example{Elementary Observation 6.1} Each path in $K\Gamma\setminus I$
is a mast of a uniserial left
$\la$-module. In other words, a sequence $\SS= (S(0),\dots,S(l))$ of simple
left $\la$-modules occurs as composition series of a uniserial module if and
only if $K\Gamma\setminus I$ contains a path through the corresponding
sequence $(e(0),\dots,e(l))$ of vertices in $\Gamma$.
\qed\endexample

\proclaim{Lemma 6.2} Suppose that $e(0),e(i),e(j)$ are vertices in $\Gamma$,
that
$\alpha : e(i)\rightarrow e(j)$ is an arrow, and that $q\in K\Gamma\setminus
I$ is a path of the form $q=q'u : e(0)\rightarrow e(j)$, where the right
subpath
$u : e(0)\rightarrow e(i)$ of
$q$ may have length zero. Moreover, assume that $\alpha u$ belongs to
$K\Gamma\setminus I$. 

If $(\alpha,u)$ is a detour on
$q$, then $(\alpha,u)$ is a slack halyard on $q$. More strongly, for each
$k\in K$, there exists a uniserial left $\la$-module $U_k$ with mast $q$
having a top element $x_k$ such that $\alpha ux_k= kqx_k$, and the
graph of $U_k$ relative to $x_k$ is
\ignore{
$$
\xy\xymatrix{
e(0) \ar@{.}[d] \ar@^{|-|}[d]<-4ex>_u \ar@^{|-|}[dd]<-8ex>_q\\
e(i) \ar@{.}[d] \ar@{-}@/^1pc/[d]^\alpha \ar@^{|-|}[d]<-4ex>_{q'}\\
e(j)}\endxy$$
}
\noindent in case $k\ne0$.\endproclaim

\demo{Proof} We will systematically forego the residue notation for
`paths in $\la$', i.e., for residue classes $p+I\in\la$, where $p\in
K\Gamma\setminus I$ is a path; it should be clear from the context whether
we are working in $K\Gamma$ or in $\la$.

Fix $k\in K$, and let $R_k\subseteq \la e(0)$ be the left ideal generated by
$\alpha u-kq$ and by the elements $r$, where $r=re(0)$ is a path in $K\Gamma$
different from $\alpha u$ and not a right subpath of $q$. Set $U_k=\la
e(0)/R_k$ and $x_k= e(0)+R_k$. To ensure the desired properties of $U_k$ and
$x_k$, it clearly suffices to show that $qx_k$ is nonzero in $U_k$. By
construction, the equality $qx_k=0$ would entail the existence of paths
$a_s$,
$b_{tr}$ in $K\Gamma$ and of scalars $k_s,l_{tr}\in K$ such that
$$q- \biggl[ \sum_s k_sa_s(\alpha u-kq) +\sum_{t,r} l_{tr}b_{tr}r\biggr] \in
I;$$
here the $a_s$ may be chosen distinct, and $r$ runs through those paths
starting in $e(0)$ which are different from $\alpha u$ and not right subpaths
of $q$. By hypothesis, $q\notin I$ and $\alpha u$ is not a right subpath of
$q$. Since $I\subseteq K\Gamma$ is generated by paths, we therefore conclude
that $q$ needs to cancel against a term $k_ska_sq$ with $k_s\ne 0$ and
$a_s=e(j)$ of length zero. But then $k_s\alpha u$ does not cancel against
any of the remaining terms inside the square brackets, since the paths $r$
are not right subpaths of $\alpha u$ either; indeed, $r\ne \alpha u$ and
each proper right subpath of $\alpha u$ is a right subpath of $q$. We infer
$\alpha u\in I$, which contradicts our hypothesis and thus yields $qx_k\ne
0$.
\qed\enddemo

Next, we will prove that the sufficient condition for finite uniserial type
given in Theorem 5.7, namely that 
$(\bold{Suf},\,p)$ be valid for all paths $p\in K\Gamma\setminus I$, is also
necessary (note that condition (N) follows from the global validity of
$(\bold{Suf},\,p)$ for all $p\in K\Gamma\setminus I$). One of the most
interesting aspects of the following theorem lies in the strong restrictions
which finiteness of the uniserial type of $\la$ imposes on the `ground plans'
for the composition series of uniserial $\la$-modules (see (4) and (4')
below). In rough terms, condition (4') says that, whenever a composition
factor of a uniserial module occurs with multiplicity larger than 1, a whole
segment of the composition series needs to repeat.

\proclaim{Theorem 6.3} For any monomial relation algebra $\la$, the following
statements are equivalent and left-right symmetric:

{\rm (1)} $\la$ has finite uniserial type.

{\rm (2)} For each sequence $\SS$ of simple left $\la$-modules, there is at
most one uniserial left $\la$-module with composition series $\SS$, up to
isomorphism.

{\rm (3)} Each mast $p$ of positive length of a uniserial left $\la$-module
satisfies $(\bold{Suf},\,p)$.

{\rm (4)} Condition (N) holds, and for each path $p\in K\Gamma\setminus I$
through a sequence of vertices $(e(0),\dots,e(l))$ in $\Gamma$, the
following is true: If $e(i)=e(j)$ and $e(i+1)\ne e(j+1)$ for some $i<j<l$,
but $K\Gamma\setminus I$ contains a path through $(e(0),\dots,e(i),e(j+1))$,
the sequence $(e(0),\dots,e(i))$ is a (left and) \underbar{right} segment of
$(e(0),\dots,e(i),\dots,e(j))$.\endproclaim

\remark{Remark} We give a coordinate-free rendering of condition (4). Namely,
condition (N) holds, and for each composition series
$\SS= (S(0),\dots,S(l))$ of a uniserial left $\la$-module $U$ the following is
true: Whenever $S(i)\cong S(j)$ and $S(i+1)\not\cong S(j+1)$ for some
$i<j<l$ such that there is a uniserial module with composition series
$(S(0),\dots,S(i),S(j+1))$, the sequence
$(S(0),\dots,S(i))$ is a (left and) \underbar{right} segment of
$(S(0),\dots,S(i),\dots, S(j))$.\endremark

\demo{Proof} Once we have shown the equivalence of (1)--(4'), left-right
symmetry of all these conditions will follow, since (1) is clearly left-right
symmetric. In view of Observations 6.1  and Theorem 5.7, we
clearly have
`(4')$\Longleftrightarrow$(4)$\Longrightarrow$(3)$\Longrightarrow$(1)'.
Condition (4) even implies that the hypothesis of Corollary 5.8 is satisfied
for each path $p\in K\Gamma\setminus I$, whence `(4)$\Longrightarrow $(2)' by
that corollary. Moreover, `(2)$\Longrightarrow$(1)' is trivial.

`(1)$\Longrightarrow$(4)'. Suppose that (1) -- and thus also condition
(N) -- holds, and let $p$ be a path through the sequence of vertices
$(e(0),\dots,e(l))$ such that $e(i)=e(j)$ for some $i<j<l$ while $e(i+1)\ne
e(j+1)$ and $K\Gamma\setminus I$ contains a path $q$ through
$(e(0),\dots,e(i),e(j+1))$. Denote the arrow from $e(i)=e(j)$ to $e(j+1)$ by
$\alpha$, that from $e(i)$ to $e(i+1)$ by $\beta$, and let $u$ and $v$ be the
paths through $(e(0),\dots,e(i))$ and
$(e(0),\dots,e(j))$, respectively. Then $q=\alpha v$, the path $\alpha
u$ belongs to $K\Gamma\setminus I$, and $(\alpha,u)$ is a detour on $q$, since
$\alpha\ne\beta$. Now Lemma 6.2 guarantees that, for each
$k\in K$, there exists a uniserial $U_k\in \lamod$, together with a top
element $x_k$, such that $\alpha ux_k= kqx_k$ and the graph of $U_k$
relative to $x_k$ is
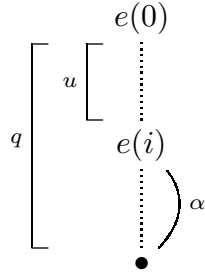
\ignore{
$$\xy\xymatrix{
e(0) \ar@{.}[d] \ar@^{|-|}[d]<-4ex>_u \ar@^{|-|}[dd]<-8ex>_q\\
e(i) \ar@{.}[d] \ar@{-}@/^1pc/[d]^\alpha\\
\bullet}\endxy$$
}
\noindent in case $k\ne 0$. In other words, if $k\ne 0$, then $(\alpha,u)$
is the only detour $(\gamma,w)$ on $q$ such that $\gamma wx_k\ne 0$. By (1),
there exist distinct scalars $k,l\in K$, with $k\ne 0$ say, such that
$U_k\cong U_l$, and [6, Section 4] yields right subpaths $w_j :
e(0)\rightarrow e(0)$ of $q$ of positive length such that
$$\alpha u \biggl( x_k+\sum_j c_jw_jx_k \biggr) =lq \biggl( x_k+\sum_j
c_jw_jx_k \biggr),$$
or, equivalently,
$ (l-k)qx_k= \sum_j c_j\alpha uw_jx_k$
since $\len(qw_j)> \len(q)$ for all $j$. The element $qx_k$ of $U_k$ being
nonzero, we derive that $\alpha uw_jx_k= k'qx_k$ for some $j$ and some
nonzero $k'\in K$. Hence the graph of $U_k$ tells us that either $\alpha uw_j
=q$, or else
$\alpha u$ is a right subpath of $\alpha uw_j$. In the latter case, we have
$\alpha uw_j= w'\alpha u$ for some path $w'$ of positive length -- because
$\len(w_j)>0$ -- and therefore $k'qx_k =\alpha uw_jx_k= w'\alpha ux_k=
kw'qx_k=0$. This being impossible, the second case is ruled out, and we
conclude that $\alpha uw_j= q= \alpha v$ in $K\Gamma$. Consequently,
$uw_j=v$, that is, $u$ is a left subpath of $v$ as required. \qed\enddemo

\head 7. Algebras of low Loewy length\endhead

Again suppose that $K$ is an infinite field. Using the results and methods of
the previous sections, one can actually list the possible shapes of uniserial
$\la$-modules of lengths bounded by a fixed integer $b$ over an algebra $\la$
which has finite uniserial type. We will provide such lists for $b=5$ and
$b=6$. In particular, they will  complete the characterization of those
algebras of finite uniserial type whose Loewy lengths are bounded above by 5
(resp., 6).

The underlying ideas are as follows: By Theorem 3.3(IId), halyards
emerging from the top vertex of a mast do not affect the isomorphism type of
the pertinent uniserial modules in the presence of condition (N). As for slack
halyards not emerging from the top vertex of a mast $p$: If $\la$ has finite
uniserial type, the presence of such halyards forces the path $p$ to run
through sequences of nested oriented cycles by Theorem 5.4.

We begin with an immediate consequence of Theorem 3.3.

\proclaim{Proposition 7.1} If $J^3=0$, then $\la$ has finite uniserial type if
and only if $\la$ satisfies condition (N). In that case, the graphs of the
uniserial $\la$-modules are either edge paths or of the form
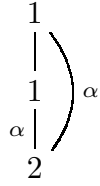
\ignore{
$$\xymatrixrowsep{1pc}
\xy\xymatrix{
1 \ar@{-}[d] \ar@{-}@/^1pc/[dd]^\alpha\\
1 \ar@{-}[d]_\alpha\\
2}\endxy$$
}
\noindent where 1 and 2 are distinct vertices. In particular, there is a 1--1
correspondence between the uniserial left $\la$-modules and their masts.
\qed\endproclaim

\proclaim{Theorem 7.2} Suppose that $\la$ has finite uniserial type, and let
$p$ be a mast of a uniserial left $\la$-module having composition length
$\le5$. Then the graphs of slack halyards on $p$ which do not emerge from the
top vertex are among the subgraphs of the following:
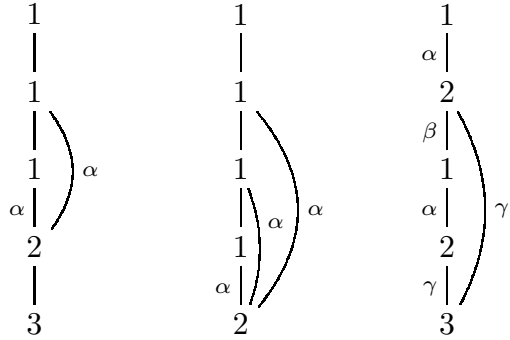
\ignore{
$$\xymatrixrowsep{1pc}
\xy\xymatrix{
1 \ar@{-}[d] &&1 \ar@{-}[d] &&1
\ar@{-}[d]_\alpha\\
1 \ar@{-}[d]
\ar@{-}@/^1pc/[dd]^\alpha &&1 \ar@{-}[d]
\ar@{-}@/^1.5pc/[ddd]^\alpha  &&2 \ar@{-}[d]_\beta
\ar@{-}@/^1pc/[ddd]^\gamma\\ 
1 \ar@{-}[d]_\alpha &&1 \ar@{-}[d]
\ar@{-}@/^0.5pc/[dd]^(0.35)\alpha &&1 \ar@{-}[d]_\alpha\\
2 \ar@{-}[d] &&1 \ar@{-}[d]_\alpha &&2 \ar@{-}[d]_\gamma\\
3 &&2 &&3}\endxy$$
}
\noindent where 1,2,3 are distinct vertices.\endproclaim

\proclaim{Corollary 7.3} For an algebra $\la$ with $J^5=0$, the following
statements are equivalent:

{\rm (1)} $\la$ has finite uniserial type.

{\rm (2)} $\la$ satisfies condition (N) of Theorem 3.3, and all slack
halyards of uniserial left $\la$-modules which do not emerge from the top
vertex of the pertinent mast must have graphs as listed in Theorem
7.2.

Moreover, if $\la$ has finite uniserial type, all of
the varieties $V_p$ are linear, and there is a bijection between
masts and isomorphism classes of uniserial left $\la$-modules.\endproclaim

\demo{Proof of Corollary 7.3} Assume (2). As mentioned above, slack halyards
emerging from the top vertex of a mast do not affect the isomorphism type of
the corresponding uniserial modules in case $\la$ satisfies (N); for
details, consult the argument establishing Theorem
3.3(IId). Moreover, that the first three types of halyards depicted in Theorem
7.2 are innocuous with respect to finiteness of the uniserial type as
well, follows from Theorem 5.1. To see that the last of the slack
halyards shown does not give rise to infinitely many uniserial modules
either, let $p= \gamma\alpha\beta\alpha$ be a mast and $(\gamma,\alpha)$ a
halyard on $p$. Moreover, let $U$ be a uniserial left $\la$-module with mast
$p$ and top element $x$ such that $\gamma\alpha x=kpx$ for some scalar $k$.
Then $y=x-k\beta\alpha x$ is in turn a top element of $U$, and $\gamma\alpha
y=0$. Since condition (N) forces the halyards to be circular, this shows
that, up to isomorphism, there exists at most one uniserial left
$\la$-module on a mast $p=\gamma\alpha\beta\alpha$ as shown in the fourth
graph of the theorem.

The other implication is an immediate consequence of the theorem. 

For a proof of the supplement, let $\la$ be of finite uniserial type, and
$p$ a mast of a uniserial left $\la$-module. Start by noting that quadratic
or higher-order terms in the polynomials defining the variety $V_p$ can only
arise if there are sequential halyards of the form
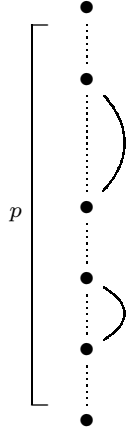
\ignore{
$$\xymatrixrowsep{1pc}
\xy\xymatrix{
\bullet \ar@{.}[d] \ar@^{|-|}[6,0]<-4ex>_p\\
\bullet \ar@{.}[dd] \ar@{-}@/^1pc/[dd]\\
\\
\bullet \ar@{.}[d]\\
\bullet \ar@{.}[d] \ar@{-}@/^1pc/[d]\\
\bullet \ar@{.}[d]\\
\bullet}\endxy$$
}
\noindent on $p$. Condition (N), however, excludes such sequences on masts
of length $\le4$ since masts of the form $1\rightarrow 1\rightarrow
2\rightarrow 2\rightarrow 3$ are disallowed. This guarantees that
$V_p$ is linear. Accordingly,
$V_p$ is isomorphic to $\AA^m$ for some $m$, where each copy of $\AA$
corresponds to a slack halyard on $p$. But relying on the list of potential
slack halyards, one easily checks that a suitable choice of top element makes
the graph of any uniserial module with mast $p$ an edge path. This shows that,
indeed, there is precisely one uniserial left $\la$-module for any given
mast.
\qed\enddemo

To complete the picture: It is not difficult to see that all the graphs of
Theorem 7.2 actually occur as graphs of slack halyards over algebras of
finite uniserial type with $J^5=0$.

\demo{Proof of Theorem 7.2} Let $p$ be a mast of length $\le4$. If the
top vertex 1 of $p$ is endowed with a loop, Theorem 5.1 tells us that
each slack halyard not emerging from the top vertex of $p$ looks as shown in
one of the three graphs. So suppose that there is no loop attached to
the vertex 1, and let $(\gamma,u)$ be a slack halyard on $p$ with
$\len(u)\ge1$. Due to the finiteness of the uniserial type of $\la$, Theorem
5.4 yields that $p$ has a right subpath $v$ of the form
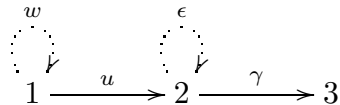
\ignore{
$$\xymatrixcolsep{3pc}
\xy\xymatrix{
1 \udotloopr{w} \ar[r]^u &2 \udotloopr{\epsilon} \ar[r]^\gamma &3
}\endxy$$
}
\noindent where 2 and 3 denote vertices, $w$ is a cycle of positive length
and $\epsilon$ a cycle which may be trivial. Since the vertex 1 does not
carry a loop, $\len(w)\ge2$, and the initial arrow $\alpha$ of $w$ ends in a
vertex different from 1. From the fact that $\len(v)\le \len(p)\le 4$, we
deduce moreover that $v=p$, that the cycle $\epsilon$ is trivial, and
$\len(w)=2$, say $w=\beta\alpha$. Moreover, we see that $u$ is an arrow, and
consequently $u=\alpha$, because $u$ is a right subpath of $p$. Thus $p$ is
of the form
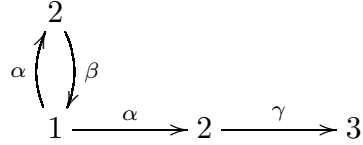
\ignore{
$$\xymatrixcolsep{3pc}
\xy\xymatrix{
2 \ar@/^0.5pc/[d]^\beta\\
1  \ar[r]^\alpha \ar@/^0.5pc/[u]^\alpha &2 \ar[r]^\gamma &3
}\endxy$$
}
\noindent The pair $(\gamma,u) =(\gamma,\alpha)$ being a detour on $p$, we
obtain $\gamma\ne\beta$; in particular, the terminal vertex 3 of $\gamma$ is
different from 1, the quiver $\Gamma$ being without double arrows. To
conclude that $3\ne2$, we note that the subpath $\alpha\beta: 2\rightarrow
1\rightarrow 2$ of $p$ is in turn a mast, whence condition (N) excludes the
existence of a loop $\gamma : 2\rightarrow 2$. Thus the vertices 1,2,3 are
distinct and the halyard $(\gamma,u)$ on $p$ is as claimed. \qed\enddemo

Without proof, we describe the slack halyards that may occur on masts $p$
of length 5 over an algebra $\la$ of finite uniserial type. They either
emerge from the top vertex of $p$ or they occur as subgraphs of the graphs
in the following two groups. The first group is
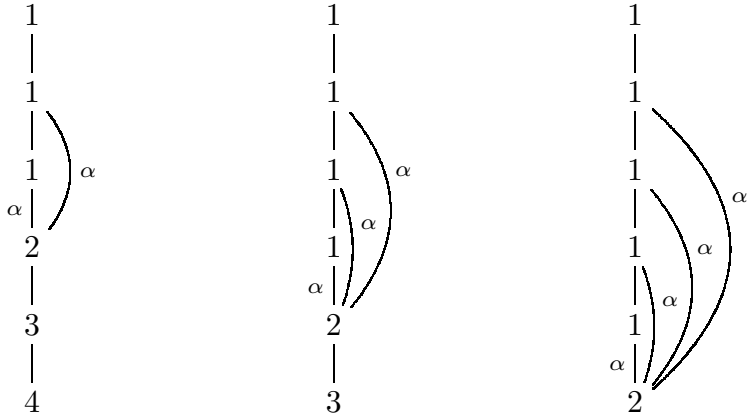
\ignore{
$$\xymatrixrowsep{1pc}
\xy\xymatrix{
1 \ar@{-}[d] &&&1 \ar@{-}[d] &&&1 \ar@{-}[d]\\
1 \ar@{-}[d] \ar@{-}@/^1pc/[dd]^\alpha &&&1 \ar@{-}[d]
\ar@{-}@/^1.5pc/[ddd]^(0.35)\alpha &&&1 \ar@{-}[d]
\ar@{-}@/^2.5pc/[dddd]^(0.35)\alpha\\
1 \ar@{-}[d]_\alpha &&&1 \ar@{-}[d]
\ar@{-}@/^0.5pc/[dd]^(0.35)\alpha &&&1 \ar@{-}[d]
\ar@{-}@/^1.5pc/[ddd]^(0.35)\alpha\\
2 \ar@{-}[d] &&&1 \ar@{-}[d]_\alpha &&&1 \ar@{-}[d]
\ar@{-}@/^0.5pc/[dd]^(0.35)\alpha\\
3 \ar@{-}[d] &&&2 \ar@{-}[d] &&&1 \ar@{-}[d]_\alpha\\
4 &&&3 &&&2}\endxy$$
}
\noindent where 1,2,3 are distinct vertices and $4\notin \{1,2\}$; these
graphs cover the case where the initial vertex 1 of $p$ is endowed with a
loop. The graphs
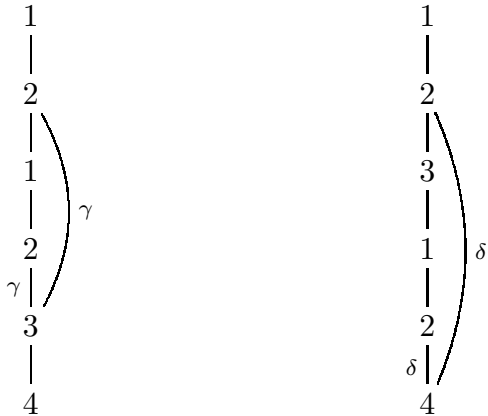
\ignore{
$$\xymatrixrowsep{1pc}
\xy\xymatrix{
1 \ar@{-}[d] &&&&1 \ar@{-}[d]\\
2 \ar@{-}[d] \ar@{-}@/^1pc/[ddd]^\gamma &&&&2 \ar@{-}[d]
\ar@{-}@/^1pc/[dddd]^\delta\\
1 \ar@{-}[d] &&&&3 \ar@{-}[d]\\
2 \ar@{-}[d]_\gamma &&&&1 \ar@{-}[d]\\
3 \ar@{-}[d] &&&&2 \ar@{-}[d]_\delta\\
4 &&&&4}\endxy$$
}
\noindent with distinct vertices 1,2,3,4 cover the case where the vertex 1
is not equipped with a loop.

For higher composition lengths, more involved patterns of slack halyards may
occur. For instance, in part (b) of Example 5.6, we presented an algebra
$\la$ of finite uniserial type which has Loewy length 7. The single
uniserial left $\la$-module with mast
$p=\alpha_5\alpha_4\alpha_3\alpha_1\alpha_2\alpha_1$ has graph
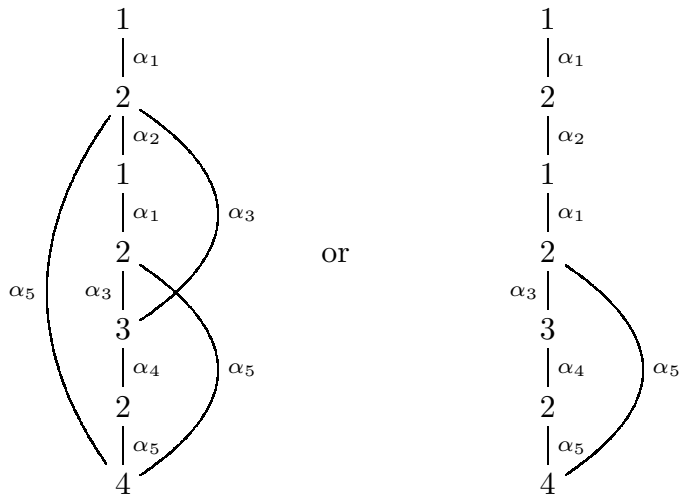
\ignore{
$$\xymatrixrowsep{1pc}
\xy\xymatrix{
1 \ar@{-}[d]^{\alpha_1} &&&&1 \ar@{-}[d]^{\alpha_1}\\
2 \ar@{-}[d]^{\alpha_2} \ar@{-}@/_2pc/[5,0]_{\alpha_5}
\ar@{-}@/^2.5pc/[3,0]^{\alpha_3} &&&&2 \ar@{-}[d]^{\alpha_2}\\
1 \ar@{-}[d]^{\alpha_1} &&&&1 \ar@{-}[d]^{\alpha_1}\\
2 \ar@{-}[d]_{\alpha_3} \ar@{-}@/^2.5pc/[3,0]^{\alpha_5}
&&\text{or}&&2 \ar@{-}[d]_{\alpha_3} \ar@{-}@/^2.5pc/[3,0]^{\alpha_5}\\
3 \ar@{-}[d]^{\alpha_4} &&&&3 \ar@{-}[d]^{\alpha_4}\\ 
2 \ar@{-}[d]^{\alpha_5}&&&&2 \ar@{-}[d]^{\alpha_5}\\ 
4 &&&&4}\endxy$$
}
\noindent according to the choice of top element; only the halyard
$(\alpha_5, \alpha_1\alpha_2\alpha_1)$ is tight.

\head Acknowledgement\endhead

The author would like to thank Axel Boldt for carefully proofreading a
preliminary manuscript of this paper.

\enddocument